\documentclass[12pt]{elsarticle}

\usepackage{amsmath, amssymb,amsthm}
\usepackage[matrix,graph]{xy}
\usepackage{comment, version}
\usepackage{color}
\newenvironment{NB}{
\color{blue}{\bf NB}.  
}{}

\excludeversion{NB}

\newtheorem{thm}{Theorem}[section]

\newtheorem{prop}[thm]{Proposition}
\newtheorem{defn}[thm]{Definition}
\newtheorem{lem}[thm]{Lemma}
\newtheorem{rem}[thm]{Remark}

\newtheorem{conj}[thm]{Conjecture}
\def\Eu{\mathop{\mathrm{Eu}}\nolimits}

\def\Spec{\mathop{\mathrm{Spec}}\nolimits}
\def\End{\mathop{\mathrm{End}}\nolimits}

\def\adj{\mathop{\mathrm{adj}}\nolimits}
\def\fund{\mathop{\mathrm{fund}}\nolimits}

\def\pt{\text{pt}}

\def\Proj{\text{Proj}}

\def\C{\mathbb{C}}
\def\R{\mathbb{R}}

\def\e{\varepsilon}
\def\mbi#1{\boldsymbol{#1}}

\newcommand{\Res}{\operatornamewithlimits{Res}}

\newcommand{\mk}{\mathfrak}
\newcommand{\mc}{\mathcal}
\newcommand{\mb}{\mathbb}

\newcommand{\mo}{\mathcal{O}}
\newcommand{\E}{\mathcal{E}}

\newcommand{\PP}{\mathbb{P}}

\newcommand{\GL}{\operatorname{GL}}

\newcommand{\wt}{\widetilde}
\newcommand{\Dec}{\operatornamewithlimits{Dec}}

\def\min{\mathop{\mathrm{min}}\nolimits}

\def\Hom{\mathop{\mathrm{Hom}}\nolimits}

\def\im{\mathop{\mathrm{im}}\nolimits}

\def\ker{\mathop{\mathrm{ker}}\nolimits}

\def\dim{\mathop{\mathrm{dim}}\nolimits}

\def\rk{\mathop{\mathrm{rk}}\nolimits}

\def\id{\mathop{\mathrm{id}}\nolimits}

\def\Q{\mathbb{Q}}

\title{Wall-crossing for vortex partition function and handsaw quiver varierty}

\author{Ryo Ohkawa \fnref{label1} \corref{cor1}}
\fntext[label1]{Osaka Metropolitan University, 
3-3-138 Sumiyoshi-ku, Osaka, 558-8585, Japan: ohkawa.ryo@omu.ac.jp, 
Research Institute for Mathematical Sciences, Kyoto University, 
Oiwake-cho, Sakyo-ku, Kyoto, 606-8502, Japan: ohkawa@kurims.kyoto-u.ac.jp}

\author{Yutaka Yoshida \fnref{label2}}
\fntext[label2]{Department of Current Legal Studies, 
Faculty of Law, Meiji gakuin university, 
1518 Kamikurata-cho, Totsuka-Ku, Yokohama
244-8539, Japan: yutakayy@law.meijigakuin.ac.jp}

\cortext[cor1]{Corresponding author}

\begin{document}

\begin{abstract}
We investigate vortex partition functions defined from integrals over the handsaw quiver varieties of type $A_{1}$ via wall-crossing phenomena.
We consider vortex partition functions defined by two types of cohomology classes, and get functional equations for each of them. 
We also give explicit formula for these partition functions. 
This gives proofs to formula suggested by physicists.
In particular, we obtain geometric interpretation of formulas for multiple hypergeometric functions including 
rational limit of the Kajihara transformation formula. 
\end{abstract}


\maketitle


\section{Introduction}
We investigate vortex partition functions defined from integrals over the handsaw quiver varieties of type $A_{1}$
via wall-crossing phenomena.
We get functional equations for two types of partition functions defined from different cohomology classes. 
We also give explicit formula for these partition functions. 
This gives proofs to formula suggested by physicists Gomis-Floch \cite{GF} and Honda-Okuda \cite{HO}
( see Section \ref{subsec:comparison}).
Furthermore, we obtain geometric interpretations of formula for multiple hypergeometric functions.
They are considered as rational limits of the Kajihara transformation \cite{K} and
formula obtained in Langer-Schlosser-Warnaar \cite{LSW} 
and Halln\"as-Langmann-Noumi-Rosengren \cite{HLNR1}, \cite{HLNR2}
from various contexts
( see Section \ref{subsec:Kajihara}).

These partition functions are analogue of the Nekrasov partition function.
The Nekrasov partition function \cite{Nek} is given by integrals over the quiver varieties of Jordan type.
These varieties are non-compact, and hence integrations are defined by counting torus fixed points.
Furthermore descriptions of fixed points depend on stability conditions for quiver representations.
From this observation, Nakajima-Yoshioka \cite{NY2} derived {\it the blow-up formula} using 
wall-crossing formula developed by Mochizuki \cite{M}.

Nekrasov's conjecture states that these partition functions give deformations of the Seiberg-Witten prepotentials for four-dimensional (4d)  $\mathcal{N}=2$ supersymmetric gauge theories.
This conjecture is proven in Braverman-Etingof \cite{BE}, Nekrasov-Okounkov \cite{NO} and Nakajima-Yoshioka \cite{NY1} independently.
In \cite{NY1}, they study relationships with similar partition functions defined for blow-up $\hat{\PP}^{2}$ of $\PP^{2}$ along the origin, and get {\it the blow-up formula}. 
These are bilinear relations of $Z_{\PP^{2}}(\e_{1}, \e_{2}-\e_{1}, \mbi a, q)$ and $Z_{\PP^{2}}(\e_{1}- \e_{1}, \e_{2}, \mbi a, q)$, which correspond to $T^{2}$-fixed points on $(-1)$ curve of $\hat{\PP}^{2}$.
Furthermore these arguments are extended in \cite{NY2} to various cohomology classes other than $1$ using the theory of perverse coherent sheaves.
This method also gives functional equations of Nekrasov functions \cite{O1}, \cite{O2}.
 

In  physics literature, vortex analogues of instanton partition functions called {\it vortex partition functions} are known. 
The vortex partition functions are given by the partition functions of the moduli data of vortex solutions in supersymmetric gauge theories, where  Higgs branch vacua  agrees with 
 the D-brane realization of the  vortex moduli space. When the gauge group is $U(r)$,  vortex moduli space in the brane construction \cite{HT}  coincides with 
handsaw quiver varieties  of  type $A_{1}$, where $r$ corresponds to the rank of a framing. 
When  supersymmetic localization is applied, one can show that vortex partition functions in two (resp. three) dimensions  have combinatorial descriptions 
similar to Nekrasov formula of  instanton partition functions in four (resp. five) dimensions. 
%
It was shown in \cite{HYY}  that, when the Fayet-Iliopoulos (FI) parameter for the moduli space  changes from  positive  to  negative,
 vortex partition functions  in certain classes of  three-dimensional (3d) $\mathcal{N}=2$ gauge theories  change discontinuously, i.e.,  
wall-crossing phenomena of Witten indices in \cite{HKY} occur.  
The authors of \cite{HYY}, see also \cite{HP}, derived   wall-crossing formulas of 
generating functions of vortex partition functions in 3d $\mathcal{N}=2$ gauge theories. 

In two-dimensional (2d) $\mathcal{N}=(2,2)$ gauge theories,
relations similar to wall-crossing formula in \cite{HYY} were originally found in the study of 
Seiberg-like dualities \cite{GF, HO}, see also earlier works \cite{BC, BSTV}.  
These relate the vortex partition functions of two different  gauge theories.  
However the precise relation between the formulas in \cite{GF, HO} and wall-crossing phenomena has not been studied  even in the physics literature.
In this paper, we study algebraic and geometric  aspects of  vortex partition functions of 2d $\mathcal{N}=(2,2)$ 
gauge theories and prove
 the wall-crossing formula of the 2d $\mathcal{N}=(2,2)$ gauge theories, which establishes the  precise relation between wall-crossing phenomena and 
 the formulas in \cite{GF, HO}.
The wall-crossing formula we will derive is regarded as a rational degeneration of the wall-crossing formula in \cite{HYY} and identical to  a finite type version of the wall-crossing formulas in \cite{O1, O2}.


Organization of the paper is the following.
In Section 2, we introduce the handsaw quiver variety of type $A_{1}$, and state our main results Theorem \ref{main0} and Theorem \ref{main}.
In Section 3, we give a proof of Theorem \ref{main0} using combinatorial description of framed moduli on $\PP^{2}$ due to 
Nakajima-Yoshioka \cite{NY1}.
In Section 4, we briefly recall Mochizuki method to prove Theorem \ref{main}.
In Section 5, we give a proof of Theorem \ref{main}.

%
%



\section{Main results}
For a fixed integer $r>0$, we consider partitions $J=(J_{0}, J_{1})$ of $[r]=\lbrace 1, \ldots, r \rbrace$
such that $[r]=J_{0} \sqcup J_{1}$, 
and put $r_{0} =\lvert J_{0} \rvert$ and $r_{1}=\lvert J_{1} \rvert$.
In our notation, we often use a vector $\vec{r}=(r_{0}, r_{1})$ instead of $J=(J_{0}, J_{1})$, and 
consider the case where 
$J_{0} = \lbrace 1, 2, \ldots, r_{0} \rbrace$ and $J_{1}=\lbrace r_{0}+1, r_{0}+2, \ldots, r \rbrace$.
For full generality, we can simply change variables.    

\subsection{Handsaw quiver variety of type $A_{1}$}
For a vector space $W=\C \mbi w_{1} \oplus \cdots \oplus \C \mbi w_{r}$, we put
\begin{align}
\label{wbase}
W_{0} = \bigoplus_{\alpha \in J_{0}} \C \mbi w_{\alpha}, 
\quad
W_{1} = \bigoplus_{\alpha \in J_{1}} \C \mbi w_{\alpha}. 
\end{align}
For a vector space $V=\C^{n}$, we consider an affine space
\begin{align*}
\mb M(W, V) &= \mb M = \End_{\C}(V) \times \Hom_{\C}(W_{0}, V) \times \Hom_{\C}(V, W_{1}).
\end{align*}
\begin{defn}
For $\zeta \in \R$, a datum $\mc A=(B, z, w) \in \mb M$ is said to be $\zeta$-stable if any sub-space $P$ of $V$ with $B(P) \subset P$ satisfies the following two conditions:
\begin{enumerate}
\item[\textup{(1)}] If $P \subset \ker w$ and $P \neq 0$, we have $\zeta \dim P < 0$.
\item[\textup{(2)}] If $\im z \subset P$ and $P \neq V$, we have $\zeta \dim V/P > 0$. 
\end{enumerate}
\end{defn}

We put $M_{+}(\vec r, n) = [ \mb M^{+} / \GL(V) ], M_{-}(\vec r, n) =  [ \mb M^{-} / \GL(V) ]$, where 
\[
\mb M^{\pm} = \lbrace (B,z,w) \in \mb M \mid \zeta \text{-stable} \rbrace
\] 
for a parameter $\zeta$ satisfying $\pm \zeta>0$.
We have an isomorphism from $M_{-}(\vec r, n)$ to {\it Laumon space} of type $A_{1}$ by \cite{N0}.
We also consider tautological bundles $\mc V= \mu^{-1}(0)^{\pm \zeta} \times V / \GL(V)$, 
$\mc W_{0} = M_{\pm}(\vec r, n) \otimes W_{0}$, 
and $\mc W_{1} = M_{\pm}(\vec r, n) \otimes W_{1}$ on $M_{\pm}(\vec r, n)$.

For an algebraic torus $T=\C^{\ast}$, we put $\mb T=T^{2r+2}$.
We write by $(q, \mbi e, e^{\mbi m}, e^{\theta})$ an element in $\mb T=T^{1} \times T^{r} \times T^{r} \times T^{1}$,
where $\mbi e=(e_{1}, \ldots, e_{r})$ and $e^{\mbi m}=(e^{m_{1}}, \ldots, e^{m_{r}})$.
We write by a monomial $m$ of these variables the corresponding character of $\mb T$, and by $\C_{m}$ the weight space, that is, 
one-dimensional $\mb T$-representation with the eigenvalue $m$.
If we rewrite
\[
\mb M(W, V) = \mb M = \End_{\C}(V) \otimes \C_{q} \times \Hom_{\C}(W_{0}, V) \times \Hom_{\C}(V, W_{1}) \otimes \C_{q},
\] 
then we have a natural $\mb T$-action on 
$M_{\pm}(\vec r, n)=\mb M^{\pm}(W, V)/\GL(V)$ by regarding $\mbi e=(e_{1}, \ldots, e_{r})$ as an diagonal elements in $\GL(W)$ via \eqref{wbase}.
Here this action does not depend on $e^{\mbi m}, e^{\theta}$, 
but in the following we will consider vector bundles on $M_{\pm}(\vec r, n)$ 
tensored with these weights.

\subsection{Partition functions}

Putting 
\begin{align}
\label{matterbdl}
\psi_{\fund } = \Eu \left( \bigoplus_{f=1}^{r} \mc V \otimes \C_{q^{-1} e^{m_{f}}} \right),
\quad
\psi_{\adj } = \Eu (TM_{\pm}(\vec r, n) \otimes \C_{e^{\theta}}),
\end{align}
we consider two types of partition functions
\begin{align*}
Z^{J}_{\pm\fund} (\e, \mbi a, \mbi m,  p)
&=  
\sum_{n=0}^{\infty} p^{n} \int_{M_{\pm}(\vec r, n)} \psi_{\fund },
\\
Z^{J}_{\pm\adj } (\e, \mbi a, \theta,  p)
&=  
\sum_{n=0}^{\infty} p^{n} \int_{M_{\pm}(\vec r, n)} \psi_{\adj }
\end{align*}
where $\mbi a=(a_{1}, \ldots, a_{r}), \mbi m = (m_{1}, \ldots, m_{r})$.

One of main results of this paper is the following explicit presentation of  these partition functions.
\begin{thm}
\label{main0}
(a) The partition function
$Z^{J}_{+\fund} (\e, \mbi a, \mbi m,  p)$ is equal to
\begin{align}
\sum_{ \substack{\mbi k \in \mathbb{Z}_{\ge 0} ^{J_{1}} } }  
p^{\lvert \mbi k \rvert}
\prod_{\beta \in J_{1}} 
\frac{ 
(-1)^{r\lvert \mbi k \rvert} \prod_{f=1}^{r}  ( (a_{\beta} + m_{f} )/\e )_{k_{\beta}}
}
{
\prod_{\alpha \in J_{1}} ( (a_{\alpha} - a_{\beta} )/\e - k_{\beta} )_{k_{\alpha}} 
\prod_{\alpha \in J_{0}} ( (a_{\alpha} - a_{\beta} )/\e - k_{\beta}  )_{k_{\beta}}
},
\label{comb2}
\end{align}
and the partition function $Z^{J}_{-\fund} (\e, \mbi a, \mbi m,  p)
$
is equal to 
\begin{align}
\sum_{ \substack{\mbi k \in \mathbb{Z}_{\ge 0} ^{J_{0}} } }  
p^{\lvert \mbi k \rvert }
\prod_{\beta \in J_{0}} 
\frac{
(-1)^{r \lvert \mbi k \rvert} 
\prod_{f=1}^{r}  ( (a_{\beta} + m_{f} )/\e  - k_{\beta} )_{k_{\beta}}
 }{\displaystyle
 \prod_{\alpha \in J_{0}} ((- a_{\alpha} + a_{\beta})/\e  - k_{\beta}  )_{k_{\alpha}} 
\prod_{\alpha \in J_{1}} ((- a_{\alpha} + a_{\beta})/\e  - k_{\beta} )_{k_{\beta}}
},
\label{comb1}
\end{align}
where $(x)_{k}= x ( x+ 1 ) \cdots (x+(k-1) )$ is the Pochhammer symbol.\\
(b)
The partition function
$Z^{J}_{+\adj} (\e, \mbi a, \mbi m,  p)
$
is equal to
\begin{align}
\sum_{ \substack{\mbi k \in \mathbb{Z}_{\ge 0} ^{J_{1}} } }
p^{\lvert \mbi k \rvert}  
\prod_{\substack{
\beta \in J_{1} \\
\alpha \in J_{1}}} 
\frac{( (a_{\alpha} - a_{\beta} - \theta )/\e - k_{\beta} )_{k_{\alpha}}}{( (a_{\alpha} - a_{\beta})/\e - k_{\beta}  )_{k_{\alpha}}} 
\prod_{\substack{
\beta \in J_{1} \\
\alpha \in J_{0}}} 
\frac{( (a_{\alpha} - a_{\beta} - \theta )/\e - k_{\beta} )_{k_{\beta}}}{( (a_{\alpha} - a_{\beta} )/\e - k_{\beta} )_{k_{\beta}}},
\label{comb4}
\end{align}
and the partition function $Z^{J}_{-\adj} (\e, \mbi a, \mbi m,  p)$
is equal to
\begin{align}
\sum_{ \substack{\mbi k \in \mathbb{Z}_{\ge 0} ^{J_{0}} } }
p^{ \lvert \mbi k \rvert}  
\prod_{\substack{
\beta \in J_{0} \\
\alpha\in J_{0}}}  
\frac{( (a_{\beta} - a_{\alpha} - \theta ) /\e - k_{\beta} )_{k_{\alpha}}}{( (a_{\beta} - a_{\alpha} )/\e - k_{\beta} )_{k_{\alpha}}} 
\prod_{\substack{
\beta \in J_{0} \\
\alpha \in J_{1}}} 
\frac{( (a_{\beta} - a_{\alpha} - \theta )/\e - k_{\beta}  )_{k_{\beta}}}{( (a_{\beta} - a_{\alpha} )/\e - k_{\beta} )_{k_{\beta}}}.
\label{comb3}
\end{align}
\end{thm}

These explicit formula are proved in Section \ref{sec:proof} using combinatorial descriptions
of fixed points sets $M_{\pm}(\vec{r}, n)^{\mb T}$.

\begin{rem}
\label{kajiharaform}
Using the elementary identity
\[
\prod_{1\le \alpha < \beta \le m} \frac{x_{\alpha}-x_{\beta} + k_{\alpha} - k_{\beta} }{x_{\alpha}-x_{\beta}}
\prod_{\alpha, \beta =1}^{m} \frac{1}{(x_{\alpha}-x_{\beta}+ 1 )_{k_{\alpha}}}
=
\prod_{\alpha, \beta =1}^{m} \frac{(-1)^{\lvert \mbi k \rvert}}{(x_{\alpha}-x_{\beta} - k_{\beta})_{k_{\alpha}}}
\]
 as in the argument preceding \cite[Theorem 2.1]{HLNR2}, the above explicit 
 formula \eqref{comb2} and \eqref{comb1} are rewritten as follows.
 
We put $u=(-1)^{(r_{1}-1)} p$.
Then $Z^{J}_{+\fund} (\e, \mbi a, \mbi m, u  )$ is equal to
\begin{align}
\sum_{ \substack{\mbi k \in \mathbb{Z}_{\ge 0} ^{J_{1}} } }  
u^{\lvert \mbi k \rvert}
\frac{\Delta_{J_{1}}(\mbi a / \e + \mbi k )}{ \Delta_{J_{1}}(\mbi a / \e)}
\prod_{\substack{ \alpha \in J_{1} \\ \beta \in J_{1}}} 
\frac{\left( \frac{a_{\alpha}}{\e} + \frac{m_{\beta}}{\e}  \right)_{k_{\alpha}}}
{\left( \frac{a_{\alpha}}{\e} - \frac{a_{\beta}}{\e} + 1 \right)_{k_{\alpha}} }
\prod_{\substack{ \beta \in J_{1} \\ \alpha \in J_{0}} } 
\frac{\left( \frac{a_{\beta}}{\e} + \frac{m_{\alpha}}{\e}  \right)_{k_{\beta}}}
{\left( \frac{a_{\beta}}{\e} - \frac{a_{\alpha}}{\e} + 1 \right)_{k_{\beta}}},
\notag
\end{align}
and $Z^{J}_{-\fund} (\e, \mbi a, \mbi m,  u)$ is equal to
\begin{align}
\sum_{ \substack{\mbi k \in \mathbb{Z}_{\ge 0} ^{J_{0}} } }  
u^{\lvert \mbi k \rvert}
\frac{\Delta_{J_{0}}(-\mbi a / \e + \mbi k )}{ \Delta_{J_{0}}(-\mbi a / \e)}
\prod_{\substack{\beta \in J_{0}\\ \alpha \in J_{0}}} 
\frac{\left( - \frac{a_{\beta}}{\e} - \frac{m_{\alpha}}{\e} + 1 \right)_{k_{\beta}}}
{\left( - \frac{a_{\alpha}}{\e} + \frac{a_{\beta}}{\e} + 1 \right)_{k_{\alpha}} }
\prod_{\substack{ \beta \in J_{0} \\ \alpha \in J_{1}}} 
\frac{\left( - \frac{a_{\beta}}{\e} - \frac{m_{\alpha}}{\e} + 1 \right)_{k_{\beta}}}
{\left( \frac{a_{\alpha}}{\e} - \frac{a_{\beta}}{\e} + 1 \right)_{k_{\beta}}}, 
\notag
\end{align}
where $\Delta_{K}(x)=\prod_{\substack{ \alpha, \beta \in K \\\alpha < \beta}} (x_{\alpha} - x_{\beta})$ 
for any subset $K \subset [r]$.
\end{rem}

Another main results are the following wall-crossing formula.
\begin{thm}
\label{main}
(a) For $A = \displaystyle \sum_{\alpha=1}^{r} a_{\alpha}/\e + m_{\alpha}/\e $,
we have 
\begin{align}
\label{fund}
Z^{J}_{+\fund}( \e, \mbi a, \mbi m,  p) = (1 + ( - 1 )^{r_{1}} p)^{- A +r_{0}} 
Z^{J}_{-\fund}(\e, \mbi a, \mbi m,  p).
\end{align}
(b) Conjecture \ref{conj} implies the following formula
\begin{align}
\label{adj}
Z^{J}_{+\adj }( \e, \mbi a, \theta,  p) = 
(1 -   p)^{(r_{0} - r_{1} )(\theta /\e + 1)} 
Z^{J}_{-\adj }(\e, \mbi a, \theta,  p).
\end{align}
\end{thm}
Here Conjecture \ref{conj} is a combinatorial identity. 
We prove this theorem in Section 5 using Mochizuki method.
In particular, we have
\begin{align*}
Z_{+\fund }^{( \emptyset, [r])}( \e, \mbi a, \mbi m,  p) 
&= 
(1 + ( - 1 )^{r} p)^{-A},
\quad
Z_{-\fund }^{([r], \emptyset)}( \e, \mbi a, \mbi m,  p) 
= 
(1 + p)^{A-r}, \\
Z_{+ \adj }^{( [r], \emptyset )}( \e, \mbi a, \theta,  p) 
&=
Z_{- \adj }^{( \emptyset, [r])}( \e, \mbi a, \theta,  p) 
=
(1 -  p)^{-r(\theta /\e +1)}.
\end{align*}

\subsection{Comparison with physical computations}
\label{subsec:comparison}
From \eqref{comb4} and \eqref{comb3} we can see that 
\[
Z^{(J_{0}, J_{1}) }_{+\adj }(\e, \mbi a, \theta, p)=Z^{(J_{1}, J_{0}) }_{-\adj } (\e, -\mbi a, \theta, p).
\]
Then our main result \eqref{fund} gives a proof of the formula proposed by Gomis-Floch \cite[(B.1.27)]{GF} after substituting 
$r=N_{F}, r_{0}=N$, and
\begin{align}
p=(-1)^{N_f-N-1}x, \quad m_f-1 /\e= {\rm i} \tilde{m}_{f}, \quad a_{f} / \e ={\rm i} m_{f},
\label{eq:ident1}
\end{align}
where the left (resp. right) hand side of \eqref{eq:ident1} denotes the symbols  in this paper (resp. \cite{GF}).  
Our main result \eqref{adj} also gives a proof of the formula proposed by Honda-Okuda \cite[(6.3)]{HO} after substituting
$r=N_{F}, r_{0}=N$, and
\[
p=-(-1)^{N_F} e^{-t}, \quad \theta/\e= - m_{ad}, \quad a_{f} / \e = m_{f}.
\]

Here we also mention on the result  in \cite{HYY}, where 
Hwang, Yi and the  second named author studied wall-crossing formulas of 3d $\mathcal{N}=2$ and $\mathcal{N}=2^*$ gauge theories in the physics literature.
In the rational limit: $\sinh(x/2) = x/2 +\cdots$, 
vortex partition functions (4.7)  (resp. (4.10)) of  a 3d $\mathcal{N}=2$  gauge theory in \cite{HYY} agree with  
\eqref{comb1} (resp. \eqref{comb2} )  under the following identification of the parameters:
\begin{align}
 &r_0=N_c, \quad  r = N_f=N_a, \quad \mbi k=(n_1, n_2, \cdots, n_{N_{c}}), \quad a_{\alpha}  =m_a,  \nonumber \\
 &\e =-2 \gamma, \quad  m_f = -\tilde{m}_\alpha + 2 \mu -2 \gamma , \quad 0 =\kappa,
\end{align}
and  also vortex partition functions (4.51) (resp. (4.54)) of  a 3d $\mathcal{N}=2^*$  gauge theory in \cite{HYY}  agree with \eqref{comb3} (resp. \eqref{comb4})   under the following identification of the parameters:
\begin{align}
 &r_0=N_c,  \quad \mbi k=(n_1, n_2, \cdots, n_{N_{c}}), \quad a_{\alpha}  =m_a,  \quad \e =-2 \gamma, \quad \theta = -2 \mu.
 \end{align}
Then the wall-crossing formulas (4.45)  and (4.59) in \cite{HYY}  can be regarded 
as the trigonometric version of  the wall-crossing formulas  studied in this article.

\subsection{Kajihara transformation}
\label{subsec:Kajihara}
Using Remark \ref{kajiharaform} and (a) in Theorem \ref{main} together, we have
\begin{align}
&
\sum_{ \substack{\mbi k \in \mathbb{Z}_{\ge 0} ^{J_{1}} } }  
u^{\lvert \mbi k \rvert}
\frac{\Delta_{J_{1}}(\mbi a / \e + \mbi k )}{ \Delta_{J_{1}}(\mbi a / \e)}
\prod_{\substack{ \alpha \in J_{1} \\ \beta \in J_{1}}} 
\frac{\left( \frac{m_{\beta}}{\e} + \frac{a_{\alpha}}{\e} \right)_{k_{\alpha}}}
{\left( 1+ \frac{a_{\alpha}}{\e} - \frac{a_{\beta}}{\e}  \right)_{k_{\alpha}} }
\prod_{\substack{ \alpha \in J_{1} \\ \beta \in J_{0}} } 
\frac{\left( \frac{m_{\beta}}{\e} + \frac{a_{\alpha}}{\e} \right)_{k_{\alpha}}}
{\left( 1+\frac{a_{\alpha}}{\e} - \frac{a_{\beta}}{\e}  \right)_{k_{\alpha}}}
\label{fundexplicit}
\\
&= 
(1 -   u)^{-A+r_{0} }
\notag
\\
&\cdot
\sum_{ \substack{\mbi k \in \mathbb{Z}_{\ge 0} ^{J_{0}} } }  
u^{\lvert \mbi k \rvert}
\frac{\Delta_{J_{0}}(- \mbi a / \e +\mbi k )}{ \Delta_{J_{0}}(-\mbi a / \e)}
\prod_{\substack{ \alpha \in J_{0}\\ \beta \in J_{0}}} 
\frac{\left( - \frac{m_{\beta}}{\e} + 1 - \frac{a_{\alpha}}{\e} \right)_{k_{\alpha}}}
{\left( 1 - \frac{a_{\alpha}}{\e} + \frac{a_{\beta}}{\e} \right)_{k_{\alpha}} }
\prod_{\substack{ \beta \in J_{0} \\ \alpha \in J_{1}}} 
\frac{\left( - \frac{m_{\alpha}}{\e} + 1 - \frac{a_{\beta}}{\e} \right)_{k_{\beta}}}
{\left( \frac{a_{\alpha}}{\e} + 1 - \frac{a_{\beta}}{\e} \right)_{k_{\beta}}}, 
\notag
\end{align}
where $u=(-1)^{(r_{1}-1)} p$, and 
$A = \displaystyle \sum_{\alpha=1}^{r} a_{\alpha}/\e + m_{\alpha}/\e $.

This formula \eqref{fundexplicit} is equivalent to the rational limit \cite[(4.8)]{K} of the Kajihara transformation
after substituting  
\[
\begin{cases}
x_{\alpha}=a_{\alpha}/ \e +1 & \alpha \in J_{1}\\ 
y_{\beta}=  - a_{\beta}/ \e & \beta \in J_{0}  
\end{cases}
\text{ and }
\begin{cases}
t_{\alpha}= -a_{\alpha}/ \e  - m_{\alpha} / \e  &  \alpha \in J_{1} \\ 
s_{\beta}= a_{\beta}/ \e + m_{\beta} / \e -1 & \beta \in J_{0}.
\end{cases}
\]
It is a natural question whether we can obtain similar geometric interpretation for the Kajihara transformation \cite[Theorem 1.1]{K}, 
or the Kajihara-Noumi transformation \cite[Theorem 2.2]{KN} in elliptic case.
See also \cite[Remark 2.4]{KN}. 

When the radius of the circle $S^1$ goes to zero,   
vortex partition functions  in 3d $\mathcal{N}=2$ theories on $S^1 \times \mathbb{R}^2$ reduce to 
vortex partition functions in 2d $\mathcal{N}=(2,2)$ theories on $\mathbb{R}^2$. In this limit,  
the wall-crossing formula  of 3d vortex partition functions \cite{HYY} reduces to that of 2d vortex partition functions 
\cite{GF} which is identical to \eqref{fund} and \eqref{fundexplicit}. 
At this moment, we have not identified 
the wall-crossing formula of the 3d vortex partition functions in \cite{HYY} with 
the original Kajihara transformation  \cite[Theorem 1.1]{K}.  
It would be interesting to clarify the relation between the wall-crossing formula of 3d vortex partition functions and   
Kajihara transformation. 

On the other hand, 
using the explicit formulas \eqref{comb4} and \eqref{comb3}, 
Theorem \ref{main} (b)
is essentially a rational version of a transformation 
formula for multiple elliptic hypergeometric series proposed in Langer-Schlosser-Warnaar
\cite[Cor. 4.3]{LSW} in the context of Kawanaka's conjecture, and Halln\"as-Langmann-Noumi-Rosengren \cite[(6.7)]{HLNR2}
in relation to deformed elliptic Ruijsenaars models.
In fact, from the special case of \cite[(6.7)]{HLNR2}, we can derive
the following formula
\begin{align}
&
\sum_{ \substack{\mbi k \in \mathbb{Z}_{\ge 0} ^{r_{1}} } }
p^{\lvert \mbi k \rvert }
\prod_{\alpha, \beta=1}^{r_{1}}  
\frac{ \displaystyle
(q^{- k_{\beta} +1 } x_{\alpha} / tx_{\beta}; q)_{k_{\alpha}}}
{ \displaystyle
(q^{-k_{\beta}} x_{\alpha} /  x_{\beta}; q)_{k_{\alpha}}}
\prod_{\beta=1}^{r_{1}}  
\prod_{\alpha=1}^{r_{0}}  
\frac{(x_{\beta} y_{\alpha}; q)_{k_{\beta}}}
{(qx_{\beta} y_{\alpha} / t; q)_{k_{\beta} }}
\label{noumi1}
\\
&=
\prod_{s=1}^{r_{1}-r_{0}} 
\frac{(q^{s}p/t^{s-1} ;q)_{\infty}}{(q^{s } p /t^{s};q)_{\infty}}
\notag\\
&\cdot
\sum_{ \substack{\mbi k \in \mathbb{Z}_{\ge 0}^{r_{0}} } }
\left( \frac{q^{r_{1}-r_{0}} p}{t^{r_{1}-r_{0}}}  \right)^{\lvert \mbi k \rvert}
\prod_{\alpha, \beta=1}^{r_{0}}  
\frac{ \displaystyle
(q^{-k_{\beta}+1} y_{\alpha} / t y_{\beta}; q)_{k_{\alpha}}}
{ \displaystyle
(q^{-k_{\beta}} y_{\alpha} / y_{\beta}; q)_{ k_{\alpha}}}
\prod_{\beta=1}^{r_{0}}  
\prod_{\alpha=1}^{r_{1}}  
\frac{(y_{\beta} x_{\alpha}; q)_{k_{\beta}}}
{(qy_{\beta} x_{\alpha} / t; q)_{k_{\beta} }}
\notag
\end{align}
due to M. Noumi (private communication).
Substituting $t=q^{\theta/\e+1}$, and
\[
\begin{cases}
x_{\alpha}=q^{a_{\alpha}/\e} & \alpha \in J_{0}\\ 
y_{\alpha}=  q^{-a_{\alpha}/\e + \theta/\e+1} & \alpha \in J_{1}
\end{cases}
\]
and taking limit $q \to 1$ in \eqref{noumi1}, we get \eqref{adj}.

\section{Proof of explicit formula}
\label{sec:proof}
Here we give proofs of explicit formulas for $Z^{J}_{\pm \fund}(\e, \mbi a, \mbi m, p)$ and $Z^{J}_{\pm \adj}(\e, \mbi a, \theta, p)$
using combinatorial description due to \cite{NY1} and \cite{N0}.
\subsection{Combinatorial description}
\label{subsec:comb}
As in \cite[Section 2 (ii)]{N0}, we can embed handsaw quiver varieties into certain framed moduli of sheaves on $\PP^{2}$ whose fixed points sets
are described by Young diagrams \cite[Proposition 2.9]{NY1}.
As a result, we also get combinatorial description of handsaw quiver varieties \cite[Section 4]{N0}.
Here we summarize it for the handsaw quiver variety of type $A_{1}$, that is, $M_{-}(\vec{r},n)$, and its
dual space $M_{+}(\vec{r},n)$.

The fixed points set $M_{-}(\vec r, n)^{\mb T}$ 
can be identified with 
$
\lbrace
\mbi k \in (\mathbb{Z}_{\ge 0})^{J_{0}} \mid \lvert \mbi k \rvert=n  
\rbrace
$,
where we regard $(\mathbb{Z}_{\ge 0})^{J_{0}}$ as a subset of $(\mathbb{Z}_{\ge 0})^{r}$ and put 
$\lvert \mbi k \rvert=\sum_{\alpha =1}^{r} k_{\alpha}$.
The $\mb T$-fixed point in $M_{-}(\vec r, n)$ corresponding to $\mbi k$ is described by
\begin{align}
\label{taut-}
V=\bigoplus_{\alpha \in J_{0}} \bigoplus_{i=1}^{k_{\alpha}} \C_{e_{\alpha} q^{-i+1}}, 
\quad
W_{0}=\bigoplus_{\alpha \in J_{0}} \C_{e_{\alpha}}, 
\quad
W_{1}=\bigoplus_{\alpha \in J_{1}} \C_{e_{\alpha}},
\end{align}
\[
B(\C_{e_{\alpha} q^{-i+1}}) =
\begin{cases}
\C_{e_{\alpha} q^{-i}} & \text{ if } 1 \le i < k_{\alpha} \\
0 & \text{ if } i = k_{\alpha} \\
\end{cases},
\quad
z(\C_{e_{\alpha}} ) =
\C_{e_{\alpha}},
\]
and $w=0$.
Hence the tangent space $T_{\mbi k}M_{-}(\vec r, n)$ at the point corresponding to $\mbi k$ is
\begin{align}
\notag
&
\Hom_{\C}( V, V) \otimes \C_{q} + \Hom_{\C}(W_{0}, V) + \Hom_{\C}(V, W_{1}) \otimes \C_{q} - \Hom_{\C}(V, V)\\
\label{tan+}
&=
\sum_{\alpha, \beta \in J_{0}} e_{\alpha} e_{\beta}^{-1} 
\left( \sum_{ i =1}^{k_{\alpha}} \sum_{j=1}^{k_{\beta}} ( q^{j -i + 1} - q^{j-i} ) 
+ \sum_{i =1}^{k_{\alpha}}  q^{-i+1} \right)
+ \sum_{\alpha \in J_{1}} \sum_{\beta \in J_{0}} e_{\alpha} e_{\beta}^{-1} \sum_{i=1}^{k_{\beta}} q^{i}.
\end{align}
Here we calculate in $K_{\mb T}(\pt)=\mathbb{Z}[q, \mbi e, e^{\mbi m}, e^{\theta}]$.

Similarly, we can identify $M_{+}(\vec r, n)^{\mb T}$ with 
$
\lbrace
\mbi k \in (\mathbb{Z}_{\ge 0})^{J_{1}} \mid \lvert \mbi k \rvert=n
\rbrace.
$
The $\mb T$-fixed point in $M_{+}(\vec r, n)$ corresponding to $\mbi k$ is described by
\begin{align}
\label{taut+}
V=\bigoplus_{\alpha \in J_{1}} \bigoplus_{i=1}^{k_{\alpha}} \C_{e_{\alpha} q^{i}}, 
\quad
W_{0}=\bigoplus_{\alpha \in J_{0}} \C_{e_{\alpha}}, 
\quad
W_{1}=\bigoplus_{\alpha \in J_{1}} \C_{e_{\alpha}},
\end{align}
\[
B(\C_{e_{\alpha} q^{i}}) =
\begin{cases}
\C_{e_{\alpha} q^{i-1}} & \text{ if } i > 1 \\
0 & \text{ if } i = 1 \\
\end{cases},
\quad
w(\C_{e_{\alpha} q^{i}}) =
\begin{cases}
0 & \text{ if } i > 1 \\
\C_{e_{\alpha}} & \text{ if } i = 1 \\
\end{cases},
\]
and $z=0$.
Hence the tangent space $T_{\mbi k} M_{+}(\vec r, n)$ at the point corresponding to $\mbi k$ is
\begin{align}
\notag
&
\Hom_{\C}(  V, V) \otimes \C_{q} + \Hom_{\C}(W_{0}, V) + \Hom_{\C}( V, W_{1}) \otimes \C_{q} - \Hom_{\C}(V, V)\\
\label{tan-}
&=
\sum_{\alpha, \beta \in J_{1}} e_{\alpha} e_{\beta}^{-1} 
\left( \sum_{ i =1}^{k_{\alpha}} \sum_{j=1}^{k_{\beta}} ( q^{i-j +1} - q^{i-j} ) 
+ \sum_{i =1}^{k_{\beta}}  q^{-i+1} \right)
+ \sum_{\alpha \in J_{1}} \sum_{\beta \in J_{0}} e_{\alpha} e_{\beta}^{-1} \sum_{i=1}^{k_{\alpha}} q^{i}.
\end{align}
Hence we have 
\begin{align}
\label{stvsco}
T_{\mbi k} M_{+}((r_{0}, r_{1}), n) = 
T_{\mbi k} M_{-}((r_{1}, r_{0}), n)) 
\Big\vert_{e_{1}=e_{1}^{-1}, \ldots, e_{r}=e_{r}^{-1}}.
\end{align}

\subsection{Proof of Theorem \ref{main0}}

Using combinatorial descriptions in the previous subsection, we rewrite $TM_{\pm}(\vec{r}, n)$ as follows.
\begin{rem}
\label{tangent1}
We have
\begin{align}
T_{\mbi k}M_{+}(\vec r, n) 
&= 
\sum_{\alpha, \beta \in J_{1}} e_{\alpha}^{-1} e_{\beta}  
\sum_{\ell=k_{\beta} - k_{\alpha}+1}^{k_{\beta}}  q^{\ell}
+
\sum_{\alpha \in J_{0}} \sum_{\beta \in J_{1}} e_{\alpha}^{-1} e_{\beta} 
\sum_{i=1}^{k_{\beta}} q^{i}.
\label{tm+}
\\
T_{\mbi k}M_{-}(\vec r, n)  
&= 
\sum_{\alpha, \beta \in J_{0}} e_{\alpha} e_{\beta}^{-1} \sum_{\ell=k_{\beta} - k_{\alpha}+1}^{k_{\beta}}  q^{\ell}
+
\sum_{\alpha \in J_{1}} \sum_{\beta \in J_{0}} e_{\alpha} e_{\beta}^{-1} 
\sum_{i=1}^{k_{\beta}} q^{i} 
\label{tm-}
\end{align}
\end{rem}
\proof 
For the first term in $T_{\mbi k}M_{-}(\vec r, n)$, we divide the range of the summation into small pieces according to $\ell =j- i=-m$.  
When $k_{\beta} \ge k_{\alpha}$, we have
\begin{align*}
&
\sum_{ i =1}^{k_{\alpha}} \sum_{j=1}^{k_{\beta}} ( q^{j -i + 1} -  q^{j-i} ) 
+ \sum_{i =1}^{k_{\alpha}} q^{-i+1} 
\\
=&
q^{k_{\beta}} + \sum_{\ell=k_{\beta} - k_{\alpha}+1}^{k_{\beta}-1} \left( \sum_{i=1}^{k_{\beta} - \ell +1 } q^{\ell} - \sum_{i=1}^{k_{\beta} - \ell} q^{\ell} \right) 
+
\sum_{\ell=1}^{k_{\beta} - k_{\alpha}} \left( \sum_{i=1}^{k_{\alpha}} q^{\ell} - \sum_{i=1}^{k_{\alpha}} q^{\ell} \right) \\
&+
\sum_{m=0}^{k_{\alpha}-1} \left( \sum_{i=m+2}^{k_{\alpha} } q^{-m} - \sum_{i=m+1}^{k_{\alpha}} q^{-m} \right) + \sum_{i =1}^{k_{\alpha}} q^{-i+1}
=
\sum_{\ell=k_{\beta} - k_{\alpha}+1}^{k_{\beta}} q^{\ell}.
\end{align*}
When $k_{\beta} < k_{\alpha}$, the similar arguments hold.
Hence in both cases, we have 
\[
e_{\alpha} e_{\beta}^{-1} \sum_{ i =1}^{k_{\alpha}} \sum_{j=1}^{k_{\beta}} 
( q^{j -i + 1} - q^{j-i} )
+ \sum_{i =1}^{k_{\alpha}} e_{\alpha} e_{\beta}^{-1} q^{-i+1} 
= e_{\alpha} e_{\beta}^{-1} \sum_{\ell=k_{\beta} - k_{\alpha}+1}^{k_{\beta}}  q^{\ell}.
\]
For $T_{\mbi k}M_{+}(\vec r, n)$, the assertion follows from \eqref{stvsco}.
\endproof

Our integral is defined by
\[
\int_{M_{\pm} (\vec{r}, n)^{\mb T} } \psi =
\sum_{\mbi k \in M_{\pm} (\vec{r}, n)^{\mb T} } \frac{\psi\vert_{\mbi k}}{\Eu( T_{\mbi k} M_{\pm}(\vec{r},n))},
\]
in the similar way as in \cite[Section 4]{NY1} ( see also Section \ref{subsec:local} for geometric background ),
where $M_{\pm} (\vec{r}, n)^{\mb T}$ are identified as subsets of $\mathbb{Z}_{\ge 0}^{r}$ as in the preceding argument.
Here the Euler class is defined by a group homomorphism $K(\text{pt}) =\mathbb{Z}[q, \mbi e, \mbi \mu] \to \mb Q (\e, \mbi a, \mbi m)^{\times}$ via 
\[
\Eu( \pm q^{i} e_{1}^{j_{1}} \cdots e_{r}^{j_{r}} \mu_{1}^{k_{1}} \cdots \mu_{r}^{k_{r}} ) =
\left(  i \e + \sum_{\alpha=1}^{r} ( j_{\alpha} a_{\alpha} + k_{\alpha} m_{\alpha} ) \right)^{\pm1}.
\]

Since combinatorial descriptions of tautological bundles $\mc V= \mu^{-1}(0)^{\pm \zeta} \times V / \GL(V)$
over each point $\mbi k \in M_{\pm}(\vec{r}, n)$
are given by \eqref{taut-} and \eqref{taut+},
this gives a proof of Theorem \ref{main0}.

\section{Mochizuki method}
Here we briefly recall methods developed in \cite{M}, and \cite{NY2} for quiver setting.
These are similar to \cite{O1}, \cite{O2}, and we will also summarize in \cite{O4} in more general setting.

\subsection{$\ell$-stability}
For any subset $\mk I$ of $[n]=\lbrace 1, \ldots, n\rbrace$ and $\C^{n}=\C \mbi v_{1} \oplus \cdots \oplus \C \mbi v_{n}$, 
we put $\C^{\mk I} = \bigoplus_{i \in \mk I} \C \mbi v_{i}$.
In this section, we fix $\mk I \subset [n]$, and often put $V=\C^{\mk I}$.
We enhance data $\mc A=(B,z,w) \in \mb M(W, \C^{\mk I})$ 
with flags $F_{\bullet} =(F_{k})$ of $\C^{\mk I}$ such
that $F_{0}=0, F_{n}= \C^{\mk I}$, and $\dim F_{k} / F_{k-1} \le 1$ for any $k$.
Hence $F_{\bullet}$ may have repetitions.
For $\ell=0,1, \ldots, n$, we introduce $\ell$-stability condition for such pairs $(\mc A, F_{\bullet})$ .
\begin{defn}
\label{ell}
A pair $(\mc A, F_{\bullet})$ of $\mc A \in \mb M(W, \C^{\mk I})$ and 
a flag $F_{\bullet}$ of $V=\C^{\mk I}$ is said to be $\ell$-stable if any sub-space $P$ of $V=\C^{\mk I}$ 
with $B(P) \subset P$ satisfies the following two conditions:
\begin{enumerate}
\item[\textup{(1)}] If $P \subset \ker w$ and $P \neq 0$, we have $P \cap F^{\ell} =0$.
\item[\textup{(2)}] If $\im z \subset P$ and $P \neq V$, we have $F^{\ell} \not\subset P$. 
\end{enumerate}
\end{defn}

To construct moduli spaces $\wt{M}^{\ell}(\vec{r},\mk I)$ of $\ell$-stable pairs,  
we put $\mk I^{\le k} = \lbrace i \in \mk I \mid i \le k  \rbrace$ and
\begin{align*}
\wt{ \mb M} = \wt{ \mb M}(W, \C^{\mk I}) = \mb M (W, \C^{\mk I})\times \prod_{k=0}^{n} 
\Hom_{\C}(\C^{\mk I^{\le k}}, \C^{\mk I^{\le k+1}}),
\end{align*}
and $G_{\mk I}= \GL(V) \times \prod_{k=1}^{n} \GL(\C^{\mk I^{\le k}})$.
For an element $f=(f_{k})_{k=0}^{n}  \in \prod_{k=0}^{n} 
\Hom_{\C}(\C^{\mk I^{\le k}}, \C^{\mk I^{\le k+1}})$, we put 
$F_{k} =f_{n} \cdots f_{k}(\C^{\mk I^{\le k}})  
\subset \C^{\mk I}$ and $F_{\bullet}=(F_{k})_{k=0}^{n}$ when all $f_{k}$ are injective.
We write by $\wt{\mb M}^{\ell}(W, \C^{\mk I})$ a subset consisting of $(\mc A, f) \in \wt{\mb M}(W, \C^{\mk I})$ such that 
all $f_{k}$ are injective and 
$(\mc A, F_{\bullet}=(\im f_{n} \cdots f_{k})_{k=0}^{n})$
is $\ell$-stable.
\begin{defn}
\label{slope}
For $\zeta \in \Q$ and $\eta =(\eta_{k})_{k=1}^{n} \in (\Q_{>0})^{n}$, an element 
$(B, z, w, f) \in \wt{\mb M}$ is said to be $(\zeta, \eta)$-semistable if 
all $f_{k}$ are injective and any sub-space $P$ of $V=\C^{\mk I}$ with $B(P) \subset P$ satisfies the following two conditions:
\begin{enumerate}
\item[\textup{(1)}] If $P \subset \ker w$ and $P \neq 0$, we have 
\[
\frac{\zeta \dim P + \sum_{k=1}^{n} \eta_{k} \dim P \cap F_{k} }{ \dim P} \le \frac{\sum_{k=1}^{n}  \eta_{k} \dim F_{ k}}{n+1}.
\]
\item[\textup{(2)}] If $\im z \subset P$ and $P \neq V$, we have 
\[
\frac{\zeta \dim V/P + \sum_{k=1}^{n} \eta_{k} \dim ( F_{ k} / P \cap F_{ k} ) }{ \dim V/P} \ge 
\frac{\sum_{k=1}^{n}  \eta_{k} \dim F_{k}}{n+1}.
\] 
\end{enumerate}
\end{defn}
The condition of injectivity for all $f_{k}$ also follows from (1) if we remove the assumption $P \neq 0$ and suitably treat the infinity. 
We put $\wt{ \mb M}^{(\zeta, \eta)}(W, \C^{\mk I}) = \lbrace (B, z, w, f) \colon (\zeta, \eta)\text{-semistable} \rbrace$.
For $\ell=1,2, \ldots, n$, we consider the following condition
\begin{align}
\label{condb}
\zeta + \sum_{k=\ell+1}^{n} \eta_{k} \cdot \dim \C^{\mk I^{\le k}} <
\frac{\sum_{k=1}^{n} \eta_{k} \cdot \dim \C^{\mk I^{\le k}} } {n+1}
< 
\zeta + \frac{\eta_{\ell} } {n+1}.
\end{align}
\begin{rem}
We assume \eqref{condb} for a fixed $\ell$.
Then we have 
\[
\wt{ \mb M}^{(\zeta, \eta)}(W, \C^{\mk I}) = \wt{\mb M}^{\ell}(\vec{r}, \C^{\mk I}).
\]
\end{rem}
\proof
When $P \subset \ker w$, we have
\[
\frac{\zeta \dim P + \sum_{k=1}^{n} \eta_{k} \dim P \cap \C^{\mk I^{\le k}}}{ \dim P}
\begin{cases}
> \zeta + \frac{\eta_{\ell}}{n+1} & \text{ if }P \cap \C^{\mk I^{\le \ell}} \neq 0 \\
< \zeta + \sum_{k=\ell+1}^{n} \eta_{k} \cdot \dim \C^{\mk I^{\le k}} 
& \text{ if }P \cap \C^{\mk I^{\le \ell}} = 0.
\end{cases}
\]
Hence in this case, the inequality in Definition \ref{slope} (1) is equivalent to $S_{0} \cap \tilde{V}^{\ell} = 0$
by \eqref{condb}.\\
\indent When $P \supset \im z$, we have
\begin{multline*}
\frac{\zeta \dim V/P + \sum_{k=1}^{n} \eta_{k} \dim \C^{\mk I^{\le k}} /P \cap \C^{\mk I^{\le k}}}{ \dim V/P}
\\
\begin{cases}
> \zeta + \frac{\eta_{\ell}}{n+1} & \text{ if }P \cap V_{I^{ \le \ell}} 
\neq V_{I^{ \le \ell}} \\
< \zeta + \sum_{k=\ell+1}^{n} \eta_{k} \cdot \dim \C^{\mk I^{\le k}} 
& \text{ if }P \cap V_{I^{ \le \ell}} = V_{I^{ \le \ell}}.
\end{cases}
\end{multline*}
Hence in this case, the inequality in Definition \ref{slope} (2) is equivalent to $P \cap V_{I^{ \le \ell}} \neq V_{I^{ \le \ell}}$
by \eqref{condb}.
\\
\indent In both cases, we get strict inequalities if $(B, z, w, f)$ is $(\zeta, \eta)$ semi-stable.
Hence we get the desired isomorphism.
\endproof
Thus we get moduli spaces $\wt M^{\ell}(\vec{r}, \mk I)$ of $\ell$-stable pairs as 
$\wt{ \mb M}^{(\zeta, \eta)}(W, \C^{\mk I}) / G_{\mk I}$.
Isomorphism classes of $\wt M^{\ell}(\vec{r}, \mk I)$ do not depend on $\mk I$, but $\dim \C^{\mk I}=
\lvert \mk I \rvert$. 
Hence we also use $\wt M^{\ell}(\vec{r}, \lvert \mk I \rvert)$ to denote this moduli space when we emphasize the number 
$\lvert \mk I \rvert$ of elements in $\mk I$, but not $\mk I$ itself.

\subsection{Enhanced master space}

\indent We consider the following two conditions on $\eta$:
\begin{align}
  \sum_{k=1}^{n} \eta_k l_{k}  \neq 0 \quad
  \text{for any $(l_{1},\dots, l_{n} )\in \mathbb{Z}^{n}\setminus \{0\}$
    with $\lvert I_{k} \rvert\le 2n^2$}
\label{2stab}\\
\eta_{m} > \dim V \sum_{k=m+1}^{n} \eta_{k} \dim F^{k}  \text{ for } m= 1, 2, \ldots, n. 
\label{condc}
\end{align}
The condition \eqref{2stab} are called $2$-stability condition in \cite{M}.
\\
\indent We take $\zeta_{-} < 0 < \zeta_{+} \in \Q$ 
and $\eta \in ( \mb Q_{>0} )^{n}$ such that
$(\zeta_{+}, \eta)$ satisfies \eqref{condb} for a fixed $\ell = 1, \ldots, n$, and $\eta$ satisfies \eqref{2stab} and \eqref{condc}.
If necessary we multiply enough divisible positive integer so that we can assume $\zeta_{\pm}$ and $\eta$ are all integers.
Then we consider ample $G$-linearizations 
\begin{align*}
\mc L_{\pm}
&=
\mo_{\wt{\mb M}} \otimes \bigotimes_{k=1}^{n} \C_{(\det_{\mk I^{\le k}})^{ \theta_{k}^{\pm} }}
\end{align*}
on $\wt{\mb M}=\wt{\mb M}(W, \mk I)$, where for any character $\chi \colon G_{\mk I} \to \C^{\ast}$ we write
by $\C_{\chi}$ the weight space. 
Here for any $\mk I' \subset [n]$, we write by $\det_{\mk I'} \colon \GL(\C^{\mk I'}) \to \C^{\ast}$ the determinant, 
and
$\theta_{n}^{\pm} = \zeta^{\pm} +\eta_{n} - \sum _{k=1}^{n-1} \eta_{k} \dim \C^{\mk I^{\le k}}$
and $\theta_{k}^{\pm}=\theta_{k}=\eta_{k} (n+1)$ for $k < n$.
See \cite{O3} for this choice of $\theta^{\pm}=(\theta_{k}^{\pm} )_{k=1}^{n}$.

We put
$\hat{ \mb M}=\hat{ \mb M}(W, \mk I) = \Proj \text{Sym} (\mc L_{-} \oplus \mc L_{+})$ and consider the semi-stable
locus $\hat{\mb M}^{ss}$ with respect to $\mo_{\hat{\mb M}}(1)$.
We define an enhanced master space by $\mc M = \hat{\mb M}^{ss}/ G$.
We consider $\C^{\ast}_{\hbar}$-action on $\mc M$ defined by
\begin{align}
\label{act0}
(\mc A, F_{\bullet},  [x_{-}, x_{+}]) \mapsto (\mc A, F_{\bullet}, [e^{\hbar} x_{-}, x_{+}]).
\end{align}

\begin{rem}
We have
\[
\mc M^{\C^{\ast}_{\hbar}} = \mc M_{+} \sqcup \mc M_{-} \sqcup 
\bigsqcup_{\mk I^{\sharp} \in D^{\ell}(\mk I)} \mc M_{\mk I^{\sharp}},
\]
where $D^{\ell} (\mk I) = \lbrace \mk I^{\sharp} \subset \mk I \mid \mk I^{\sharp} \neq \emptyset, \min(\mk I^{\sharp}) \le \ell \rbrace $.
\end{rem}
\indent We see that $\mc M_{-} = \lbrace x_{+}=0 \rbrace$ is isomorphic to the full flag bundle 
$Fl(\mc V, \mk I)$ of the tautological bundle $\mc V$ on $M_{-}(\vec{r}, \lvert \mk I \rvert)$
with $\lbrace k \in [n] \mid \mc F^{k} / \mc F^{k-1} \rbrace =\mk I$, and $\mc M_{+} = \lbrace x_{-}=0 \rbrace$ is 
isomorphic to $\wt{M}^{\ell}(\vec{r}, \mk I)$.

We call $\mk I^{\sharp} \in D^{\ell}(\mk I)$ decomposition data, and we identify them with a pair $(\mk I^{\sharp}, \mk I^{\flat})$, where 
$ \mk I^{\flat} = \mk I \setminus \mk I^{\sharp}$.
We take a decomposition $\C^{\mk I}=\C^{\mk I^{\sharp}} \oplus \C^{\mbi{\mk I}^{\flat}}$, and
describe $\mc M_{\mk I^{\sharp}}$ as follows.

\subsection{Modified $\C^{\ast}_{\hbar}$-action}
We put $D=\zeta^{+} - \zeta^{-} \in \mathbb{Z}$.
For each $\mk I^{\sharp} \in \mc D^{\ell}(\mk I)$, we put $d^{\sharp}=\lvert \mk I^{\sharp} \rvert$.
To describe $\mc M_{\mk I^{\sharp}}$, we consider a modified action 
$\C^{\ast}_{\frac{\hbar}{d^{\sharp}D}} \times \mc M \to \mc M$ induced by 
\begin{align}
\label{act}
\left( \mc A, F_{\bullet}, [x_{-},x_{+}] \right) \mapsto \left( e^{\frac{h}{d^{\sharp}D}} \id_{\C^{\mk I^{\sharp}} \oplus } 
\id_{\C^{\mk I^{\flat}}} \right) \left( \mc A, F_{\bullet}, [e^{\hbar}x_{-},x_{+}] \right). 
\end{align}
This action is equal to the original $\C^{\ast}_{\hbar}$-action \eqref{act0}, since the difference is absorbed in $G$-action.
Then $(\mc A_{\sharp} \oplus \mc A_{\flat}, f_{\sharp} \oplus f_{\flat}, [1, \rho])$ is fixed 
by this $\C^{\ast}_{\hbar}$-action for $(\mc A_{\sharp}, f_{\sharp}) \in \wt{\mb M}(0, \C^{\mk I^{\sharp}})$, 
$(\mc A_{\flat}, f_{\flat}) \in \wt{\mb M}(W, \C^{\mk I^{\flat}})$, 
and $\rho \neq 0$.
Hence if this element is stable, then it represents a $\C^{\ast}_{\hbar}$ -fixed point in $\mc M$.
This observation together with analysis for stability condition in \cite{O1} implies that
$\mc M_{\mk I^{\sharp}} $ is isomorphic to the quotient stack of
\[
\wt{ \mb M}^{\sharp}(0, \C^{\mk I^{\sharp}}) \times ( \C_{t^{-1} (\det_{\mk I^{\sharp}})^{D}} )^{\times}
\times
\wt{ \mb M}^{\min(\mk I^{\sharp})-1}(W, \C^{\mk I^{\flat}}) \times ( \C_{t(\det_{\mk I^{\flat}})^{D}})^{\times} 
\] 
by a group $G_{\mk I^{\sharp}} \times G_{\mk I_{\infty}} \times \C^{\ast}_{t}$,  
where $\wt{\mb M}^{\sharp} (0, \C^{\mk I^{\sharp}}) =
\lbrace (B, F_{\bullet}) \in \wt{\mb M} (0, \C^{\mk I^{\sharp}}) \mid \C[B] F^{1} = \C^{\mk I^{\sharp}} \rbrace$, and the superscript
``$\times$'' denotes the complement of zero.
\\
\indent If we replace the group $\C^{\ast}_{t}$ with $\C^{\ast}_{t^{1/d^{\sharp}D}}$, then we have an \'etale cover 
$\Phi_{\mk I^{\sharp}} \colon \mc M'_{\mk I^{\sharp}} \to \mc M_{\mk I^{\sharp}}$ of degree $1/d^{\sharp}D$. 
Then by a group automorphism of $\GL(V_{\mk I^{\sharp}}) \times \C^{\ast}_{t^{1/d^{\sharp}D}}$ 
sending $(g, t^{1/d^{\sharp}D})$ to $(t^{-1/d^{\sharp}D}g, t^{1/d^{\sharp}D})$, we get an isomorphism
$\mc M'_{\mk I^{\sharp}} \cong \mc M_{\sharp} \times \mc M_{\flat}$, where
\begin{align}
\mc M_{\sharp} &= [ \wt{\mb M}^{\sharp}(0, \C^{\mk I^{\sharp}}) \times ( \C_{(\det_{\mk I^{\sharp}})^{D}} )^{\times} / G_{\mk I^{\sharp}} ],
\label{msharp}
\\
\mc M_{\flat} &= [\wt{\mb M}^{\min(\mk I^{\sharp})-1} (W, \C^{\mk I^{\flat}}) \times ( \C_{t(\det_{\mk I^{\flat}})^{D}} )^{\times}
/ G_{\mk I^{\flat}} \times \C^{\ast}_{t^{1/d^{\sharp}D}}].
\label{mflat}
\end{align}

\subsection{Moduli stack of destabilizing objects $\mc M_{\sharp}$}
To study $\mc M_{\sharp}$, we give explicit description of $M_{-}((1, 0), d^{\sharp})=\mb M^{-}(W_{\sharp}, V_{\sharp})/\GL(V_{\sharp})$, where
$W_{\sharp}$ is a $\mathbb{Z}/2 \mathbb{Z}$-graded vector space with $W_{\sharp 0 }=\C$ and $W_{\sharp 1} = 0$.
In fact, we do not need the next proposition, but we give a proof for the sake of explanation.

\begin{prop}
\label{destab}
We have an isomorphism $M_{-}((1, 0), d^{\sharp} ) \cong \mb A^{d^{\sharp}}=\Spec \C[x_{1}, \ldots, x_{d^{\sharp}}]$ such that the tautological bundle $\mc V_{\sharp}$ corresponds to $\C[x_{1}, \ldots, x_{d^{\sharp}}]$-module $\C[x_{1}, \ldots, x_{d^{\sharp}}, y] / (y^{d^{\sharp}} + x_{1} y^{d^{\sharp}-1} + \cdots + x_{d^{\sharp}})$ via this isomorphism.
Furthermore $q \in T$ acts by $q^{-1} y$ and $q^{i} x_{i}$ for $i=1,\ldots, d^{\sharp}$.
\end{prop}
\proof
We regard $\C[x_{1}, \ldots, x_{d^{\sharp}}]$-module $\C[x_{1}, \ldots, x_{d^{\sharp}}, y] / (y^{d^{\sharp}} + x_{1} y^{d^{\sharp}-1} + \cdots + x_{d^{\sharp}})$ as family of data with $(B, z, w)=(y, 1, 0)$.
Then we have a morphism $\mb A^{d^{\sharp}}=\Spec \C[x_{1}, \ldots, x_{d^{\sharp}}] \to M_{-}((1, 0), d^{\sharp})$, and we can check that this induces bijection between the sets of closed points.
Hence this is an isomorphism.
 
We consider the corresponding torus action on $\mb A^{d^{\sharp}}$ as follows. 
For $\mbi x=(x_{1}, \ldots, x_{d^{\sharp}}) \in \mb A^{d^{\sharp}}$, we put $f_{\mbi x} = y^{d^{\sharp}} + x_{1} y^{d^{\sharp}-1} + \cdots + x_{d^{\sharp}}$.
We have a commutative diagram:
\[
\xymatrix{
\C[y] / \langle f_{\mbi x}(y) \rangle \ar[r]^{ty} & \C[y] / \langle f_{\mbi x}(y) \rangle\\
\C[y] / \langle f_{\mbi x}(t^{-1}y) \rangle \ar[u]^{y=ty} \ar[r]^{y} & \C[y] / \langle f_{\mbi x}( t^{-1}y)  \rangle \ar[u]^{y=ty} }
\]
If we put $q . \mbi x=(q x_{1}, q^{2} x_{2}, \ldots, q^{d^{\sharp}-1} x_{d^{\sharp}-1}, q^{d^{\sharp}} x_{d^{\sharp}})$, 
we have $f_{\mbi x}(q^{-1}y) = q^{-d^{\sharp}} f_{q . \mbi x} (y)$.
Hence we get the assertion for $T$-action on $\mb A^{d^{\sharp}}$.

Similarly we can easily see that $GL(W_{\sharp})$ trivially acts on $M_{-}((1, 0) , d^{\sharp})$ by the following commutative diagram
for $e \in \GL(W_{\sharp})$:
\[
\xymatrix{
\C \ar[r]^-{e^{-1}} &\C[y] / \langle f_{\mbi x}(y) \rangle \ar[r]^{y} & \C[y] / \langle f_{\mbi x}(y) \rangle\\
\C \ar@{=}[u] \ar[r]^-{1} & \C[y] / \langle e f_{\mbi x}(y) \rangle \ar[u]^{e^{-1}} \ar[r]^{y} & \C[y] / \langle e f_{\mbi x}( y)  \rangle \ar[u]^{e^{-1}} }
\]
\endproof

In the following, we write by $\mc V_{\sharp}$ the tautological bundle over $M_{-}((1, 0), d^{\sharp})$, and by $\mc V_{\flat}$ 
one on $\wt{M}^{\min(\mk I^{\sharp})-1}(\vec{r}, \mk I^{\flat})$.
We consider the full flag bundle 
$Fl(\mc V_{\sharp}/ \mo_{M_{-}((1, 0), d^{\sharp})}, \bar{\mk I}^{\sharp})$ of
the tautological bundle $\mc V_{\sharp}$ on $M_{-}((1, 0), d^{\sharp})$,
where $\bar{\mk I}^{\sharp}= \mk I^{\sharp} \setminus \lbrace \min(\mk I^{\sharp}) \rbrace$.
We also write by the same letter $\mc V_{\sharp}$ the pull-back to $Fl(\mc V_{\sharp} / \mo_{M_{-}((1, 0), d^{\sharp})}, \bar{\mk I}^{\sharp})$
and consider the determinant line bundle $\det \mc V_{\sharp}$ over $Fl(\mc V_{\sharp} / \mo_{M_{-}((1, 0), d^{\sharp})}, \bar{\mk I}^{\sharp})$.
Then for $\mc M_{\sharp}$ defined in \eqref{msharp}, we have the following proposition.
\begin{rem}
\label{dest}
We have an isomorphism
\[
\mc M_{\sharp} \cong 
[((\det \mc V_{\sharp})^{\otimes D} \otimes \C_{u^{d^{\sharp}D}} )^{\times} / \C^{\ast}_{u}].
\]
In particular, we have an \'etale morphism $\mc M_{\sharp} \to Fl(\mc V_{\sharp} / \mo_{M_{-}((1, 0), d)}, \bar{\mk I}^{\sharp})$
of degree $1/(d^{\sharp}D)$.
\end{rem}
\proof
We use \eqref{msharp} and consider a map 
\[
\wt{ \mb M}^{\sharp}(0, V_{\mk I^{\sharp}}) \to \mb M^{-}(W_{\sharp}, V_{\sharp}),
\quad
(B_{\sharp}, f_{\sharp}=(f_{\sharp k}) ) \mapsto
(B_{\sharp}, f_{\sharp n-1} \cdots f_{\sharp \min(\mk I^{\sharp})}, 0).
\]
Then a group automorphism $\GL(V_{\mk I^{\sharp}}) \times \GL(V_{\mk I^{\sharp}}^{\min(\mk I^{\sharp})})$ 
mapping $(g, u)$ to $(u^{-1}g, u)$ gives morphism $\mc M_{\sharp} \to M^{-}((1,0), d)$, which induce the desired isomorphism.
\endproof

We consider the universal flags $\bar{\mc F}_{\bullet}^{\sharp}$ on 
$Fl(\mc V_{\sharp} / \mo_{M_{-}((1, 0), d)}, \bar{\mk I}^{\sharp})$ 
and $\mc F_{\bullet}^{\flat}$ on $\wt{ M}^{\min(\mk I^{\sharp})-1}(\vec r, \mk I^{\flat})$.
On $Fl(\mc V_{\sharp} / \mo_{M_{-}((1, 0), d)}, \bar{\mk I}^{\sharp})$,
we define a full flag $\mc F_{\bullet}^{\sharp}$ of $\mc V^{\sharp}$ 
by the pull-back of $\bar{\mc F}_{i}^{\sharp}$ to $\mc V^{\sharp}$ for $i \neq \min(\mk I^{\sharp})$ 
and $\mc F_{\min(\mk I^{\sharp})}^{\sharp} = \im (\mo \to \mc V^{\sharp})$. 
We also write their pull-backs to the product 
$Fl( \mc V_{\sharp} / \mo_{M_{-}((1, 0), d)}, \bar{\mk I}^{\sharp}) \times 
\wt M^{\min(\mk I^{\sharp})-1}(\vec{r}, \mk I^{\flat} )$
by the same letter $\mc F^{\sharp}_{\bullet}, \mc F^{\flat}_{\bullet}$.

By \eqref{msharp} and \eqref{mflat} and Proposition \ref{dest}, we have a natural \'etale morphism 
$\Psi_{\mk I^{\sharp}} \colon \mc M'_{\mk I^{\sharp}} \to 
Fl(\mc V_{\sharp} / \mo_{M_{-}((1, 0), d)}, \bar{\mk I}^{\sharp} ) \times 
\wt{ M}^{\min(\mk I^{\sharp})-1}(\vec r, \mk I^{\flat})$ 
of degree $1/(d^{\sharp}D)^{2}$.

\subsection{Decomposition of $\mc M^{\C^{\ast}_{\hbar}}$}
Summarizing, we have the following theorem.
\begin{thm}
\label{decomp}
We have
\[
\mc M^{\C^{\ast}_{\hbar}} = \mc M_{+} \sqcup \mc M_{-} \sqcup 
\bigsqcup_{\mk I^{\sharp} \in \mc D^{\ell}(\mk I)}\mc M_{\mk I^{\sharp}}
\]
such that the followings hold.\\
(i) We have $\mc M_{+} \cong \wt{M}^{\ell}(\vec r, \mk I)$ and $\mc M_{-} \cong \wt{M}^{0}(\vec r, \mk I)$, that is, 
the full flag bundle $Fl(\mc V, \mk I)$ of the tautological bundle $\mc V$ over $M_{-}(\vec r, \lvert \mk I \rvert)$.\\
(ii) For each $\mk I^{\sharp} \in \mc D^{\ell}$, we have finite \'etale morphisms 
$
\Phi_{\mk I^{\sharp}} \colon \mc M'_{\mk I^{\sharp}} \to \mc M_{\mk I^{\sharp}}$ of degree $1/(d^{\sharp} D)$, and 
$\Psi_{\mk I^{\sharp}} \colon \mc M_{\mk I^{\sharp}}' \to  
Fl( \mc V_{\sharp} / \mo_{M_{-}((1, 0), d^{\sharp})}, \bar{\mk I}^{\sharp}) \times 
\wt M^{\min(\mk I^{\sharp})-1}(\vec{r}, \mk I^{\flat} )
$ 
of degree $1/(d^{\sharp}D)^{2}$, 
where $d^{\sharp}=\lvert \mk I^{\sharp} \rvert$ and $D= \zeta_{+} - \zeta_{-} \in \mathbb{Z}$. \\
(iii) As $\C^{\ast}_{\hbar/d^{\sharp}D}$-equivariant vector bundles on $\mc M_{\mk I^{\sharp}}'$, 
we have 
\[
\Phi_{\mk I^{\sharp}}^{\ast} \mc V_{\mk I^{\sharp}} \cong \Psi_{\mk I^{\sharp}}^{\ast} \mc V_{\sharp} \otimes e^{\frac{\hbar}{d^{\sharp}D}} \otimes 
(L_{\mk I^{\sharp}})^{\vee}, \quad 
\Phi_{\mk I^{\sharp}}^{\ast} \mc V_{\mk I_{\infty}} \cong \Psi_{\mk I^{\sharp}}^{\ast} \mc V_{\flat}
\]
for a line bundle $L_{\mk I^{\sharp}}$ on $\mc M_{\mk I^{\sharp}}'$ such that 
$
(L_{\mk I^{\sharp}})^{\otimes d^{\sharp}D} \cong 
\Psi_{\mk I^{\sharp}}^{\ast}( \det \mc V_{\sharp } \otimes \det \mc V_{\flat })^{\otimes D}$. 
\end{thm}

We write by $N_{+}, N_{-}$, and $N_{\mk I^{\sharp}}$ normal bundles of $\mc M_{+}, \mc M_{-}$ and $\mc M_{\mk I^{\sharp}}$ in $\mc M$ respectively.
For the following computations, we introduce some notation and describe these normal bundles here.
For a vector bundle $\mc E$ and a finite set $\mk I^{\sharp} \subset \mathbb{Z}$ with $\rk \E$ elements, 
we write by $Fl(\mc E, \mk I^{\sharp})$ the full flag bundle 
$F_{\bullet}=(F_{i})_{i \in \mathbb{Z}}$ of $\mc E$ such that 
\[
\lbrace i \in \mathbb{Z} \mid F_{i}/F_{i-1} \neq 0 \rbrace = \mk I^{\sharp}.
\]
These are all isomorphic to $Fl(\mc E, [\rk \mc E])$, but we use products of these flag bundles
and combinatorial description in Section \ref{subsec:euler2}.

For two flags $\mc F_{\bullet}, \mc F_{\bullet}'$ of sheaves on the same Deligne-Mumford stack,  
we put 
\begin{align*}
\Theta (\mc F_{\bullet}, \mc F_{\bullet}') 
&=
\sum_{i > j }  
\mc Hom\left(\mc F_{j} / \mc F_{j-1}, \mc F_{i}' /\mc F_{i-1}' \right),
\\
\wt{\mk H} (\mc F_{\bullet}, \mc F_{\bullet}')
&=
\Theta (\mc F_{\bullet}, \mc F_{\bullet}') 
+
\Theta (\mc F_{\bullet}', \mc F_{\bullet}).
\end{align*}
When $\mc F_{\bullet} = \mc F_{\bullet}'$, we put $\Theta(\mc F_{\bullet}) = \Theta(\mc F_{\bullet}, \mc F_{\bullet})$.

\begin{lem} 
\label{lem:normal}
(1) We have $N_{\pm} = \mc L_{\pm} \otimes e^{\hbar}$.\\
(2) We have $\Phi_{\mk I^{\sharp}}^{\ast} N_{\mk I^{\sharp} } \cong 
\Psi_{\mk I^{\sharp}}^{\ast} \left( \mathfrak{N}( \mc V_{\sharp} 
\otimes e^{\frac{\hbar}{d^{\sharp}D}} \otimes 
(L_{\mk I^{\sharp}})^{\vee}, \mc V_{\flat}) + 
\wt{\mk H}
( \mc F^{\sharp}_{\bullet} \otimes e^{\frac{\hbar}{d^{\sharp}D}}, \mc F^{\flat}_{\bullet}) \right)
$ where
\begin{align}
\label{normalbdl}
\mathfrak{N}( \mc V_{\sharp}, \mc V_{\flat})
&=
\mc H om( \mc V_{\flat}, \mc V_{\sharp} \otimes \C_{q}) + \mc H om( \mc W_{0}, \mc V_{\sharp}) \\  
&+ 
\mc H om( \mc V_{\sharp}, \mc V_{\flat} \otimes \C_{q}) + \mc H om( \mc V_{\sharp}, \mc W_{1} \otimes \C_{q})  \notag \\
&-
\mc H om( \mc V_{\flat}, \mc V_{\sharp}) - \mc H om( \mc V_{\sharp}, \mc V_{\flat}) \notag
\end{align}
for vector bundles $\mc V, \mc V_{\sharp}$, and $\mathbb{Z}_{2}$-graded vector bundle $\mc W=\mc W_{0} \oplus \mc W_{1}$.
\end{lem}

\subsection{Localization}
\label{subsec:local}
We put $M_{0} = \Spec \Gamma(\mb M(W, V), \mo_{\mb M})^{\GL(V)}$.
By the main result \cite[(1)]{GP} and Theorem \ref{decomp}, we have the following diagram
\[
\xymatrix{
\varprojlim_{m} A^{\ast}_{\C^{\ast}_{\hbar}} (\mc M \times_{\mb T} E_{m}) \otimes_{\C[\hbar]} \C[\hbar^{\pm 1}] \ar[d]_{\Pi_{\ast} (\cdot) \cap {[\mc M]^{vir}}} \ar[r]^{\cong}& \varprojlim_{m} A^{\ast} (\mc M^{\C^{\ast}_{\hbar}} \times_{\mb T} E_{m}) \otimes_{\C} \C[\hbar^{\pm 1}] \ar[d]_{\Pi_{\ast} (\cdot) \cap \left( [\mc M_{+}]^{vir} +  [\mc M_{-}]^{vir} + \sum_{\mk I^{\sharp}} [\mc M_{\mk I^{\sharp}}]^{vir} \right)} \\
\varprojlim_{m} A_{\ast} (M_{0} \times_{\mb T} E_{m}) \otimes_{\C} \C[\hbar, \hbar^{-1}] \ar@{=}[r] & \varprojlim_{m} A_{\ast} (M_{0} \times_{\mb T} E_{m}) \otimes_{\C} \C[\hbar,\hbar^{-1}]\\
}
\]
where the upper horizontal arrow is given by 
\[
\frac{\iota_{+}^{\ast}}{\Eu (\mk N(\mc M_{+}))} + \frac{\iota_{-}^{\ast}}{\Eu (\mk N(\mc M_{-}))} + 
\sum_{\mk I^{\sharp}} \frac{\iota_{\mk I^{\sharp}}^{\ast}}{\Eu (\mk N(\mc M_{\mk I^{\sharp}}))}.
\]
Here $\hbar$ corresponds to the first Chern class in $A^{\C^{\ast}_{\hbar}}(\text{pt})$ of the weight $e^{\hbar} \in \C^{\ast}_{\hbar}$, and $\iota_{\pm}$ and $\iota_{\mk I^{\sharp}}$ are embeddings of $\mc M_{\pm}$ and $\mc M_{\mk I^{\sharp}}$ into $\mc M$.

Then $\Eu(N_{+}), \Eu(N_{\mk I^{\sharp}})$ are invertible in $A^{\ast}_{\tilde{\mb T} \times \C^{\ast}_{\hbar}}(\mc M)[\hbar,\hbar^{-1}]]$, 
and we have the following localization formula
\begin{align}
\label{localization}
\int_{\mc M} \varphi = \int_{\mc M_{+}} \frac{\varphi\vert_{\mc M_{+}}}{\Eu (N_{+})} + 
\int_{\mc M_{-}} \frac{\varphi\vert_{\mc M_{-}}}{\Eu (N_{-})} + 
\sum_{\mk I^{\sharp}\in \mc D^{\ell}(\mk I)} 
\int_{\mc M_{\mk I^{\sharp}}} \frac{\varphi\vert_{\mc M_{\mk I^{\sharp}}}}{\Eu (N_{\mk I^{\sharp}})}.
\end{align}

\subsection{Cohomology classes}
We take two cohomology classes
\begin{align}
\label{cls}
\Eu \left( \bigoplus_{f=1}^{r} \mc V 
\otimes \frac{e^{m_{f}}}{q} \right), \quad 
\Eu^{\theta} \left( \Lambda(\mc V) \right) 
\in A^{\bullet}_{\tilde{\mb T} \times \C^{\ast}_{\hbar}} (\mc M), 
\end{align}
where 
$
\Lambda(\mc V) = \mc End(\mc V) \otimes \C_{q} + \mc Hom (W_{0}, \mc V ) 
+ \mc Hom (\mc V, W_{1}) \otimes \C_{q}  - \mc End( \mc V).
$ 
For $\alpha = [E] - [F] \in K(\mc M)$ with vector bundles $E, F$ on $\mc M$, we define $\Eu^{\theta}(\alpha) = c_{\rk E} (E \otimes e^{\theta}) / c_{\rk F} (F \otimes e^{\theta})$.

In the following we consider one of classes in \eqref{cls} and write it by $\psi = \psi (\mc V)$.
Furthermore, we put 
\begin{align*}
\tilde{\psi} = \tilde{\psi}(\mc V) = 
\frac{ \psi \cdot \Eu^{\theta} (\Theta(\mc F_{\bullet}))}{\lvert \mk I \rvert !} 
\in A_{\C^{\ast}_{\hbar}}^{\ast}( \mc M \times_{\mb T} E_{m})  
\end{align*}
and substitute $\varphi = \tilde{\psi}$ in \eqref{localization}.
The left hand side in \eqref{localization} is a polynomial in $\hbar$ while the right hand side has a power series part in $\hbar^{-1}$.
Hence if the symbol $\displaystyle \Res_{\hbar = \infty}$ denotes the operation taking the coefficient in $\hbar^{-1}$, we have
\begin{align*}
\int_{\wt{M}^{\ell}(\vec{r}, \mk I)} \tilde{\psi}
-\int_{M_{-}(\vec{r}, \lvert \mk I \rvert)} \psi 
\notag
&=
- \Res_{\hbar=\infty} \sum_{\mk I^{\sharp} \in D^{\ell}(\mk I)} \int_{\mc M_{\mk I^{\sharp}}} 
\frac{\tilde{\psi}\vert_{\mc M_{\mk I^{\sharp}}}}{\Eu (N_{\mk I^{\sharp}})}\\
\end{align*}

By Lemma \ref{lem:normal} (2),
the last summand is equal to
\begin{multline}
\label{euler}
\frac{\lvert \mk I^{\flat} \rvert!}{\lvert \mk I \rvert!} 
\int_{\wt{M}^{\min(\mk I^{\sharp})-1} (\vec{r}, \mk I^{\flat})} 
\frac{1}{\lvert \mk I^{\flat} \rvert!}
\int_{Fl(\mc V_{\sharp} / \mo_{M_{-}((1, 0), d^{\sharp} )}, \bar{\mk I}^{\sharp})} 
\\
\frac{
\psi(\mc V_{\flat} \oplus \mc V_{\sharp} \otimes e^{\hbar}) \cdot
\Eu^{\theta}(\Theta( \mc F^{\sharp}_{\bullet} \otimes e^{\hbar} \oplus \mc F^{\flat}_{\bullet} ))
} 
{
\mk N(\mc V_{\sharp} \otimes e^{\hbar}, \mc V_{\flat}) \cdot
\Eu \left( \wt{\mk H}( \mc F^{\sharp}_{\bullet} \otimes e^{\hbar}, \mc F^{\flat}_{\bullet} ) \right)
},
\end{multline}
where $\bar{\mk I}^{\sharp} = \mk I^{\sharp} \setminus \lbrace \min{\mk I^{\sharp}} \rbrace$.
We have  
$\Theta( \mc F^{\sharp}_{\bullet} \otimes e^{\hbar} \oplus \mc F^{\flat}_{\bullet} )=
\Theta(\mc F_{\bullet}^{\flat}) + \Theta(\mc F_{\bullet}^{\sharp}) + 
\wt{\mk H}( \mc F^{\sharp}_{\bullet} \otimes e^{\hbar}, \mc F^{\flat}_{\bullet} )$,
and
\begin{align}
\label{thetabar}
\Theta(\mc F_{\bullet}^{\sharp})=\Theta(\bar{\mc F}_{\bullet}^{\sharp}) + 
\mc V_{\sharp}/ \mo_{M_{-}((1,0), d^{\sharp})}.
\end{align} 
We note that $\Theta(\bar{\mc F}_{\bullet}^{\sharp})$ is the relative tangent bundle of
$Fl(\mc V_{\sharp} / \mo_{M_{-}((1,0), d^{\sharp})}, \bar{\mk I}^{\sharp}) \to M_{-}((1,0), d^{\sharp})$.

In this expression, we deleted some line bundles and a parameter $d^{\sharp}D$, since we have 
$\displaystyle \Res_{\hbar = \infty} f(\hbar) = d^{\sharp}D \Res_{\hbar = \infty} f( d^{\sharp}D \hbar+a)$ (cf. \cite[Section 8.2]{O1}), 
and integrals are taken over $1/d^{\sharp}D$-degree \'etale covering of full flag bundles
$Fl(\mc V_{\sharp} / \mo_{M_{-}((1, 0), d^{\sharp})}, \bar{\mk I}^{\sharp})$.


\section{Wall-crossing formula}
In the following, we use wall-crossing formula deduced from analysis in the previous section to get functional equations \eqref{fund}.
These are the similar calculations to \cite[Section 6]{NY2}, hence we omit detail explanation.

For $\mbi d= (d_{1}, \ldots, d_{j}) \in \mathbb{Z}_{>0} ^{j}$, we put $\lvert \mbi d \rvert = d_{1} + \cdots + d_{j}$.
Let $\Dec_{j}^{n}$ be the set of collections $\mbi{\mk I}=( \mk I_{1}, \ldots, \mk I_{j})$ of non-empty subsets of 
$[n]$ such that
\begin{enumerate}
\item[$\bullet$] $\mk I_{i} \cap \mk I_{k} = \emptyset$ for $i \neq k$, and
\item[$\bullet$] $\min(\mk I_{1}) > \cdots > \min(\mk I_{j})$.
\end{enumerate}
We note that $\Dec_{1}^{n} = \mc D^{n}$.
We identify $\mbi{\mk I}=(\mk I_{1}, \ldots, \mk I_{j})$ with $(\mk I_{1}, \ldots, \mk I_{j}, \mk I_{\infty})$ where
$\mk I_{\infty} = [n] \setminus \bigsqcup_{i=1}^{j} \mk I_{i}$.
We consider maps $\sigma \colon \Dec_{j+1}^{n} \to \Dec_{j}^{n}$,
\[
\mbi{\mk I}'= (\mk I_{1}', \ldots, \mk I_{j+1}') \mapsto \sigma(\mbi{\mk I}') = (\mk I_{1}', \ldots, \mk I_{j}')
\]
and 
$\rho \colon \text{Dec} _{j}^{n} \to S_{j}^{n}=\left \lbrace \mbi d= (d_{1}, \ldots, d_{j}) \in \mathbb{Z}_{>0} ^{j} \ 
\Big \vert \ \lvert \mbi d \rvert \le n \right \rbrace, 
\quad
\mbi{\mk I} \mapsto \mbi d_{\mbi{\mk I}} = (\lvert \mk I_{1} \rvert, \ldots, \lvert \mk I_{j} \rvert)$ 
for $j=1, \ldots, n$.


\subsection{Iterated cohomology classes}
\label{subsec:iter1}
\indent For $j > 0 $ and $\mbi d =(d_{1}, \ldots, d_{j})\in \mathbb{Z}_{>0}^{j}$, we consider a product 
\[
M_{\mbi d} = M_{-}(\vec{r}, n - \lvert \mbi d \rvert) \times \prod_{i=1}^{j} M_{-}((1, 0), d_{i}).
\]
We write by $\mc V^{(i)}$ the tautological bundle
on the component $M_{-}((1, 0), d_{i})$ of $M_{\mbi d}$.

For $\mbi{\mk I}=( \mk I_{1}, \ldots, \mk I_{j}) \in \rho^{-1}(\mbi d)$, 
we also put
\[
\wt M_{\mbi{\mk I}}^{\ell} = 
\wt M^{\ell} (\vec{r}, \mk I_{\infty}) \times
Fl(\mc V_{\sharp}^{(1)} / \mo_{M((1,0),d_{1})}, \bar{\mk I}_{1}) 
\times \cdots \times 
Fl(\mc V_{\sharp}^{(j)} / \mo_{M((1,0), d_{j})}, \bar{\mk I}_{j})
\]
where $\bar{\mk I}_{i}=\mk I_{i} \setminus \lbrace \min(\mk I_{i})\rbrace$,
and write by $\mc F_{\bullet}^{(i)}$ the pull-back to $\wt M_{\mbi{\mk I}}^{\ell}$ of 
the universal flag on each component 
$Fl(\mc V^{(i)} / \mo_{M_{-}((1,0), d_{i})}, \bar{\mk I}_{i})$.
We consider tori $\C^{\ast}_{\hbar_{i}}$ for $i =1,2, \ldots, n$, and take a coordinate $e^{\hbar_{i}} \in \C^{\ast}_{\hbar_{i}}$.
By abuse of notation $e^{\hbar_{i}}$ also denotes a trivial bundle with $e^{\hbar_{i}}$-weight.
We write by $\mc F_{\bullet}^{\infty}$ the pull-back of
the flag $\mc F_{\bullet}$ on $\wt M^{\ell}(\vec{r}, \mk I_{\infty} )$, and
put $\mc F_{\bullet}^{> i}=\mc F_{\bullet}^{\infty} \oplus \bigoplus_{k > i} \mc F_{\bullet}^{(k)} \otimes e^{\hbar_{k}}$.

We do not need the obstruction theory to define fundamental cycle 
$[\wt M_{\mbi{\mk I}}^{\ell}]$ since we only consider smooth Deligne-Mumford stacks in this paper.
For $\alpha \in A_{\mb {T} \times \C^{\ast}_{\theta} \times \prod_{i=1}^{j} \C^{\ast}_{\hbar_{i}}}^{\bullet}(\wt M_{\mbi{\mk I}}^{\ell})$, 
we write by $\int_{[\wt M_{\mbi{\mk I}}^{\ell}]} \alpha  \in 
A_{\mb {T} \times \C^{\ast}_{\theta} \times \prod_{i=1}^{j} \C^{\ast}_{\hbar_{i}}}^{\bullet}
(\wt M^{\ell}(\vec{r}, n - \lvert \mbi d_{\mbi{\mk I}} \rvert ) )$ 
the Poincare dual of the push-forward of $\alpha \cap [\wt M_{\mbi{\mk I}}^{\ell}]$ by the projection 
$\wt M_{\mbi{\mk I}}^{\ell} \to \wt M^{\ell}(\vec{r}, \mk I_{\infty} )$. 

We write by the same letters $\mc V, \mc W_{0}, \mc W_{1}$  
the pull-backs of ones on $\wt M^{\ell}(\vec{r}, \mk I_{\infty} )$ to the product $\wt M_{\mbi{\mk I}}^{\ell}$.
For $\mbi{\mk I}=(\mk I_{1}, \ldots, \mk I_{j}) \in \text{Dec} _{j}^{n} $, we write by $\tilde{\psi}_{\mbi{\mk I}} (\mc V)$
the following cohomology class
\begin{align}
&
\int_{[\wt M_{\mbi{\mk I}}^{\ell} ]} 
\frac{ 
\psi \left( \mc V \oplus \bigoplus_{i=1}^{j} \mc V^{(i)} \otimes e^{\hbar_{i}} \right) 
\Eu^{\theta}(\Theta(\mc F_{\bullet}^{> 0}))
}
{ 
\Eu \left( \bigoplus_{i=1}^{j} \mathfrak{N} ( 
\mc V^{(i)} \otimes e^{\hbar_{i}}, 
\mc V \oplus \bigoplus_{k=i+1}^{j} 
\mc V^{(k)} \otimes e^{\hbar_{k}})
\right)
}
\frac{
1
}
{
\lvert \mk I_{\infty} \rvert!
}
\label{itcoho}
\\
&=
\int_{[\wt M_{\mbi{\mk I}}^{\ell}]} 
\frac{ 
\psi \left( \mc V \oplus \bigoplus_{i=1}^{j} \mc V^{(i)} \otimes e^{\hbar_{i}} \right) 
\Eu^{\theta}(\Theta(\mc F^{\infty}_{\bullet}))
/\lvert \mk I_{\infty} \rvert!
}
{ 
\Eu \left( 
\bigoplus_{i=1}^{j} 
\mk N (\mc V^{(i)} \otimes e^{\hbar_{i}},
\mc V \oplus \bigoplus_{k=i+1}^{j} \mc V^{(k)} \otimes e^{\hbar_{k}})\right)
}
\notag
\\
&\cdot
\prod_{i=1}^{j}
\frac{
\Eu^{\theta} ( \Theta(\mc F_{\bullet}^{(i)}) ) 
\Eu^{\theta}(\wt{\mk H} (\mc F_{\bullet}^{(i)} \otimes e^{\hbar_{i}}, \mc F_{\bullet}^{> i}))
}
{
\Eu(\wt{\mk H} (\mc F_{\bullet}^{(i)} \otimes e^{\hbar_{i}}, \mc F_{\bullet}^{> i}))
}
\notag
\end{align}
in $A^{\bullet}_{\mb{T} \times \C^{\ast}_{\theta} \times \prod_{i=1}^{j} \C^{\ast}_{\hbar_{i}}}(M)$.
Here $\mathfrak{N}( \mc V_{\sharp}, \mc V_{\flat})$ is defined by
\eqref{normalbdl}.
By modified $\C^{\ast}_{\hbar}$-action \eqref{act}, we need to multiply $\mc V^{(i)}$ with $e^{\hbar_{i}}$ in \eqref{itcoho}.

By the projection formula, we have a polynomial $f_{\mbi{\mk I} }(\mbi x) \in \mb Q(\e, \mbi a, \mbi m, \theta)[\mbi x]$ 
independent of $M$ such that $f_{\mbi{\mk I}}(\mc V) = \psi_{\mbi{\mk I} }(\mc V)$  
since we have finite $\mb T$-fixed points sets of $Fl(\mc V^{(i)} / \mo_{M_{-}((1,0),d_{i})}, \bar{\mk I}_{i})$.


\subsection{Recursions}
\label{subsec:loca}
For $\mbi{\mk I} = (\mk I_{1}, \ldots, \mk I_{j}) \in \mathbb{Z}_{>0}^{j}$, we put $\mk I = \mk I_{\infty}$, 
$\ell = \min(\mk I_{j} ) -1$, 
and take an equivariant cohomology classes $\varphi = \tilde{\psi}_{\mbi{\mk I }}(\mc V)$ on $\mc M$.
For the convenience, we also put $\tilde{\psi}_{()}=\psi$ for $j=0$. 
Then \eqref{euler} is equal to
\begin{align}
\label{exc}
\frac{\lvert \mk I^{\flat} \rvert!  }{\lvert \mk I \rvert !} \int_{\wt M^{\min(\mk I^{\sharp})-1}(\vec{r}, \mk I^{\flat} ) } 
\tilde{\psi}_{ ( \mbi{\mk I}, \mk I^{\sharp} )}(\mc V),
\end{align}
where $(\mbi{\mk I}, \mk I^{\sharp}) \in \Dec^{n}_{j + 1}$.
Using this argument repeatedly, we deduce recursion formula.
\begin{lem}
For $l \ge 1$, we have
\begin{align}
\notag
&
\int_{M_{+}(\vec r, n)} \psi(\mc V)  - \int_{M_{-}(\vec{r}, n )} \psi(\mc V) \\
\notag
&= 
\sum_{j=1}^{l-1} 
\Res_{\hbar_{1} = \infty} \cdots \Res_{\hbar_{j}=\infty}
\sum_{\mbi{\mk I} \in \Dec_{j}^{n}} \frac{ \lvert \mk I_{\infty} \rvert!  }{ n! }
\int_{\wt M^{0}(\vec{r}, \mk I_{\infty} )} \tilde{\psi}_{\mbi{\mk I}} (\mc V)\\
\label{formula2>0}
&+
\Res_{\hbar_{1} = \infty} \cdots \Res_{\hbar_{l}=\infty}
\sum_{\mbi{\mk I} \in \Dec_{l}^{n}} \frac{\lvert \mk I_{\infty} \rvert!  }{ n! }
\int_{\wt M^{\min(\mk I_{j})-1}(\vec{r}, \mk I_{\infty} )} \tilde{\psi}_{\mbi{\mk I}} (\mc V).
\end{align}
\end{lem}
\proof
We prove by induction on $j$.
For $l=1$, \eqref{formula2>0} is nothing but \eqref{exc} for $j_{0}=0$ and $\ell = n$.
For $l \ge 1$, we assume the formulas \eqref{formula2>0}.
Then again by \eqref{exc}, the last summand for each $\mbi{\mk I} \in \Dec_{l}^{n}$ is equal to 
\begin{align*}
&
 \frac{\lvert \mk I_{\infty} \rvert! }{ n! }
\left( \int_{\wt M^{0}(\vec{r}, \mk I_{\infty} )}  \tilde{\psi}_{\mbi{\mk I}} (\mc V)  \right. \\
&+
\Res_{\hbar_{l+1}=\infty}
\left. \sum_{\mbi{\mk I}' \in \sigma^{-1} (\mbi{\mk I})} 
\frac{\lvert \mk I'_{\infty} \rvert! }{\lvert \mk I_{\infty} \rvert!} 
\int_{\wt M^{\min (\mk I_{l+1}' )-1 }(\vec{r}, {\mk I}'_{\infty} )} \tilde{\psi}_{\mbi{\mk I}'} (\mc V)  \right ),
\end{align*}
where $\mbi{\mk I}'=(\mbi{\mk I}, \mk I_{l+1}')$ and $\mk I'_{\infty}= 
[n] \setminus \bigsqcup_{i=1}^{l+1} \mk I'_{i}$.
Hence we have \eqref{formula2>0} for general $l \ge 1$.
\endproof

For $l> n$, the set $\Dec_{l}^{n}$ is empty.
Thus we get the following theorem.
\begin{thm}
\label{thm:wcm}
We have
\begin{align}
\label{main>0}
&
\int_{M_{+}(\vec r, n)}  \psi(\mc V)  - \int_{M_{-}(\vec r, n)}  \psi(\mc V) 
\notag\\
&=
\sum_{j=1}^{n} \sum_{\mbi{\mk I} \in \Dec_{j}^{n}} \frac{ \lvert \mk I_{\infty} \rvert! }{ n! }
\Res_{\hbar_{1} = \infty} \cdots \Res_{\hbar_{j}=\infty}
\int_{ \wt M^{0}(\vec{r}, \mk I_{\infty} )} 
\tilde{\psi}_{\mbi{\mk I}} (\mc V).
\end{align}
\end{thm}

\subsection{Euler classes of tautological bundles}
\label{subsec:euler1}
For $\psi=\Eu( \bigoplus_{f=1}^{r} \mc V \otimes e^{m_{f}}/q)$, we substitue $\theta=0$ in \eqref{euler} and \eqref{itcoho}.
So for $\mbi d=(d_{1}, \ldots, d_{j}) \in ( \mathbb{Z}_{> 0} )^{j}$, we put 
\begin{align*}
\psi_{\mbi d} (\mc V)
&=
\int_{\prod_{i=1}^{j} M_{-}((1,0), d_{i})}
\frac{ 
\psi \left( \mc V \oplus \bigoplus_{i=1}^{j} \mc V^{(i)} \otimes e^{\hbar_{i}} \right) 
\Eu \left( \bigoplus_{i=1}^{j} \mc V^{(i)}/ \mo_{M_{\mbi d}^{\ell}}  \right) 
}
{
\Eu \left( \bigoplus_{i=1}^{j} 
\mk N (\mc V^{(i)} \otimes e^{\hbar_{i}}, \mc V \oplus \bigoplus_{k=i+1}^{j} \mc V^{(k)} \otimes e^{\hbar_{k}} )\right),
}
\end{align*}
where $\int_{\prod_{k=1}^{j} M_{-}((1,0), d_{k})}$ 
denotes the push-forward by the projection 
\[
M_{-}(\vec{r}, n -\lvert \mbi d_{\mbi{\mk I}} \rvert ) \times \prod_{i=1}^{j} M_{-}((1,0), d_{i}) \to M_{-}(\vec{r}, n -\lvert \mbi d_{\mbi{\mk I}} \rvert ).
\]

Since $\wt M^{0}(\vec{r}, \mk I_{\infty} )$ is a full flag bundle of $\mc V$ over $M_{-}(\vec{r}, n - \lvert \mbi d_{\mbi{\mk I}} \rvert )$,
we have
\begin{align}
&
\int_{M_{+}(\vec r, n)} \psi (\mc V)  - \int_{M_{-}(\vec{r}, n )} \psi (\mc V) 
\notag\\
&=
\sum_{j=1}^{n} 
\sum_{\mbi{\mk I} \in \Dec_{j}^{n}} 
\frac{ \lvert \mk I_{\infty} \rvert! \prod_{k=1}^{j} (\lvert \mk I_{k} \rvert-1)! }{ n! }
\Res_{\hbar_{1} = \infty} \cdots \Res_{\hbar_{j}=\infty}
\int_{ M_{-}(\vec{r}, n - \lvert \mbi d_{\mbi{\mk I}} \rvert )} \psi_{\mbi d_{\mbi{\mk I}}} (\mc V).
\label{formula3}
\end{align}
Each summand in the right hand side of \eqref{formula3} depends only on $\mbi d=\mbi d_{\mbi{\mk I}}$.
\begin{lem}
For $\mbi d \in S_{j}^{n}=\left \lbrace \mbi d= (d_{1}, \ldots, d_{j}) \in \mathbb{Z}_{>0} ^{j} \ \Big\vert \ \lvert \mbi d \rvert \le n \right \rbrace$ and 
$\mbi{\mk I}=(\mk I_{1}, \ldots, \mk I_{j}) \in \rho^{-1}(\mbi d)$, we have 
\begin{align*}
\lvert \rho^{-1}(\mbi d) \rvert
&=
\frac{1}{ \prod_{i=1}^{j} \sum_{1 \le k \le i} d_{k}} 
\frac{ n! }{ \lvert \mk I_{\infty} \rvert! \prod_{k=1}^{j} (\lvert \mk I_{k} \rvert-1)! }.
\end{align*}
\end{lem}
\proof
This follows from \cite[Lemma 6.8]{NY2} since $\lvert \mk I_{k} \rvert = d_{k}$.
\endproof

By this lemma and \eqref{formula3}, we get 
\begin{align}
&
\int_{M_{+}(\vec r, n)} \psi (\mc V)  - \int_{M_{-}(\vec{r}, n )} \psi (\mc V) 
\notag\\
&=
\sum_{j=1}^{n} \sum_{\substack{ \mbi d \in \mathbb{Z}_{>0}^{j} \\ \lvert \mbi d \rvert \le n}} 
\frac{1}{\prod_{i=1}^{j} \sum_{1 \le k \le i} d_{k}} 
\Res_{\hbar_{1} = \infty} \cdots \Res_{\hbar_{j}=\infty}
\int_{M_{-}(\vec{r}, n - \lvert \mbi d \rvert)} \psi_{\mbi d} (\mc V).
\label{formula3'}
\end{align}

Since $\int_{M_{-}((1,0),d_{i})}  \Eu \left(\mc V^{(i)}/ \mo_{M_{-}((1,0), d_{i})} \right) 
=\frac{(-1)^{d_{i} -1}}{d_{i} \e}$, we have
\begin{align}
\psi_{\mbi d} (\mc V)
&=
\psi ( \mc V ) 
\cdot
\prod_{i=1}^{j} 
\frac{(-1)^{d_{i} -1}}{d_{i} \e}
\frac{ \Eu \left( \sum_{f=1}^{r} \mc V^{(i)} \otimes e^{\hbar_{i}} \otimes \mb{C}_{q^{-1} e^{m_{f}}} \right) } 
{
\Eu \left( 
\mk N (\mc V^{(i)} \otimes e^{\hbar_{i}}, \mc V \oplus \bigoplus_{k=i+1}^{j} \mc V^{(k)} \otimes e^{\hbar_{k}})
\right).
}
\notag
\end{align}
By \eqref{matterbdl} and \eqref{normalbdl}, 
the residue $\Res_{\hbar_{1} = \infty} \cdots \Res_{\hbar_{j}=\infty} \psi_{\mbi d}(\mc V)$ is equal to
\begin{align}
&
\Res_{\hbar_{1} = \infty} \cdots \Res_{\hbar_{j}=\infty}
\psi ( \mc V ) 
\cdot
\prod_{i=1} ^{j} 
\frac{(-1)^{d_{i} -1}}{d_{i} \e}
\prod_{l=0}^{d_{i}-1} 
\frac{
\prod_{f=1}^{r} ( \hbar_{i} - (l+1) \e  + m_{f}) 
}
{\displaystyle
\prod_{\alpha\in J_{0}} 
(\hbar_{i} - l \e  - a_{\alpha})
\prod_{\alpha \in J_{1}} 
(-\hbar_{i} + (l+1) \e  + a_{\alpha})
}
\notag \\
&=
\psi(\mc V) \cdot 
\prod_{i=1} ^{j} 
\frac{(-1)^{r_{1}d_{i}+d_{i} -1}}{\e} \cdot
\left( -\sum_{\alpha=1}^{r} (a_{\alpha} + m_{\alpha}) + r_{0} \e\right)
\label{residue}
\end{align}
since we have $\Res_{\hbar_{1} = \infty} \cdots \Res_{\hbar_{j}=\infty}
\frac{\Eu \left( 
\mc V^{>i} \otimes e^{- \hbar_{i} + l \e} + (\mc V^{>i})^{\vee} \otimes e^{\hbar_{i} - l \e}
\right)}
{
\Eu \left( 
\mc V^{>i} \otimes e^{-\hbar_{i} + (l+1) \e} + (\mc V^{>i})^{\vee} \otimes e^{\hbar_{i} - (l-1)\e} 
\right)
}
=0$, 
where we put $\mc V^{>i}=\mc V \oplus \bigoplus_{k=i+1}^{j} \mc V^{(k)} \otimes e^{\hbar_{k}}$.
Here \eqref{residue} is equal to $\psi(\mc V) \cdot 
(-1)^{ r_{1} \lvert \mbi d \rvert + \lvert \mbi d \rvert-j} (- A + r_{0})^{j}
$ where $A=\sum_{\alpha=1}^{r} (a_{\alpha} + m_{\alpha})$. 
If we put 
$
\alpha_{n}=\int_{M_{-}(\vec r, n)} \psi$ and $\beta_{n}=\int_{M_{+}(\vec r, n)} \psi,
$ 
then by \eqref{formula3'} we have
\begin{align}
\beta_{n} 
&=
\alpha_{n}+
\sum_{j=1}^{n} \sum_{\substack{ \mbi d \in (\mathbb{Z}_{>0})^{j} \\ \lvert \mbi d \rvert \le n}} \frac{(-1)^{r_{1} \lvert \mbi d \rvert + \lvert \mbi d \rvert - j} 
( - A + r_{0})^{j}}{\prod_{i=1}^{j} \sum_{1 \le k \le i} d_{k}} \alpha_{n - \lvert \mbi d \rvert}
\notag\\
&= 
\sum_{l=0}^{n} \frac{(-1)^{r_{1} l} ( - A +r_{0}) ( - A + r_{0} -1) \cdots ( - A + r_{0} - l +1)}{l!} \alpha_{n-l}. 
\notag
\end{align}
Since $Z^{J}_{-\fund} (\e, \mbi a, \mbi m,  p)= \sum_{n=0}^{\infty} \alpha_{n} p^{n}$ and 
$Z^{J}_{+\fund} (\e, \mbi a, \mbi m,  p)= \sum_{n=0}^{\infty} \beta_{n} p^{n}$, this gives \eqref{fund}.

\subsection{Euler classes of tangent bundles}
\label{subsec:euler2}
To show \eqref{adj}, let us take a class 
\[
\Lambda(\mc V) = \mc End(\mc V) \otimes \C_{q} + \mc Hom (W_{0}, \mc V ) + \mc Hom (\mc V, W_{1}) \otimes \C_{q} - \mc End(\mc V)
\in K_{\mb T}(\mc M),
\] 
and put $\psi(\mc V)=\Eu^{\theta}(\Lambda(\mc V))$ in $A^{\bullet}_{\C^{\ast}_{\theta} \times \C^{\ast}_{\hbar}}(\mc M)$.
We also put $u= (r_{1} - r_{0}) (\theta/\e +1)$.
In \eqref{euler}, we have
\begin{align*}
\tilde{\psi}(\mc V_{\flat} \oplus \mc V_{\sharp} \otimes e^{\hbar})
&=
\tilde{\psi}(\mc V_{\flat}) \psi (\mc  V_{\sharp} \otimes e^{\hbar}) 
\Eu^{\theta} \left( S(\mc V_{\flat}, \mc V_{\sharp}) \otimes \C_{q} - 
S(\mc V_{\flat}, \mc V_{\sharp}) \right),
\\
\mk N( \mc V_{\sharp} \otimes e^{\hbar}, \mc V_{\flat})
&= 
S(\mc V_{\flat}, \mc V_{\sharp}) \otimes \C_{q}  - 
S(\mc V_{\flat}, \mc V_{\sharp}) 
+ \mc Hom(W_{0}, \mc V_{\sharp}) + \mc Hom (\mc V_{\sharp}, W_{1}) \otimes \C_{q},
\end{align*}
where $S(\mc V_{\flat}, \mc V_{\sharp}) = \mc Hom(\mc V_{\flat}, \mc V_{\sharp}) + \mc Hom(\mc V_{\sharp}, \mc V_{\flat})$. 

To compute the push-forward $ \int_{Fl(\mc V_{\sharp} / \mo_{M_{-}((1, 0), p)}, \bar{\mk I}^{\sharp})}$
in \eqref{euler},
we divide the computations into four cases:
\begin{lem}
We have the following.\\
(1) 
$
\displaystyle
\Res_{\hbar = \infty} 
\Eu^{\theta} (S(\mc V_{\flat}, \mc V_{\sharp}) \otimes \C_{q}  - S(\mc V_{\flat}, \mc V_{\sharp}))/
\Eu (S(\mc V_{\flat}, \mc V_{\sharp}) \otimes \C_{q}  - S(\mc V_{\flat}, \mc V_{\sharp}))=0.
$
\\
(2) 
$
\displaystyle
\Res_{\hbar = \infty} 
\Eu^{\theta}(\wt{\mk H}(\mc F^{\sharp}_{\bullet} \otimes e^{\hbar}, \mc F^{\flat}_{\bullet})) / 
\Eu((\wt{\mk H}(\mc F^{\sharp}_{\bullet} \otimes e^{\hbar}, \mc F^{\flat}_{\bullet}))) = - s(\mk I^{\flat}, \mk I^{\sharp}) \theta, 
$
where $\mk I^{\flat}=[n] \setminus \mk I^{\sharp}$ for $\mk I^{\sharp} \subset [n]$, and 
\[
s(\mk I^{\sharp}, \mk I^{\flat}) = 
\left \lvert \lbrace (i, j) \in \mk I^{\sharp} \times \mk I^{\flat}\mid i < j \rbrace \right\rvert - 
\left\lvert \lbrace (i, j) \in \mk I^{\sharp} \times \mk I^{\flat}  \mid i > j \rbrace \right\rvert.
\]
(3) 
$
\displaystyle 
\Res_{\hbar = \infty} \frac{\Eu^{\theta}(\mc Hom(W_{0}, \mc V_{\sharp}) +\mc Hom (\mc V_{\sharp}, W_{1}) \otimes \C_{q})}
{\Eu(\mc Hom(W_{0}, \mc V_{\sharp}) + \mc Hom (\mc V_{\sharp}, W_{1}) \otimes \C_{q})} = (r_{0}-r_{1}) d^{\sharp} \theta.
$\\
(4) 
$
\displaystyle
\int_{Fl(\mc V_{\sharp} / \mo_{M_{-}((1, 0), d^{\sharp})}, \bar{\mk I^{\sharp}})} 
\Eu^{\theta}(\mc End(\mc V_{\sharp}) \otimes \C_{q} - 
\mc End (\mc V_{\sharp}) + \Theta(\mc F^{\sharp}_{\bullet} )= \frac{(\theta/\e+1)_{d^{\sharp}}}{d^{\sharp} \theta}.
$
\end{lem}
\proof
We only give a proof for (4) since the other statements are direct computations.
To prove (4) we remark that from \eqref{thetabar} we have 
\begin{align*}
\mc End(\mc V_{\sharp}) \otimes \C_{q} - 
\mc End (\mc V_{\sharp}) + \Theta_{\sharp}
&=
TM_{-}((1,0), d^{\sharp}) - \mo_{M_{-}((1,0), d^{\sharp})} + \Theta(\bar{\mc F}^{\sharp}_{\bullet}).
\end{align*}
Then the assertion follows from 
$\int_{M_{-}((1, 0), d^{\sharp})} \Eu^{\theta}( TM_{-}((1, 0), d^{\sharp})) = \frac{(\theta/\e+1)_{d^{\sharp}}}{d^{\sharp}!}$, and 
\[
\int_{Fl(\C^{\bar{\mk I}^{\sharp}}, \bar{\mk I}^{\sharp})} \Eu^{\theta}( T Fl(\C^{\bar{\mk I}^{\sharp}}, \bar{\mk I}^{\sharp}))=(d^{\sharp}-1)!.
\]
Here the last integral does not depend on $\theta$ since $Fl(\C^{\bar{\mk I}^{\sharp}}, \bar{\mk I}^{\sharp})$ is compact.
\endproof

Computing \eqref{euler} by this proposition, we have
\begin{align}
\notag
&
\int_{\wt{M}^{\ell}(\vec r, \mk I )} \tilde{\psi} 
- \int_{M_{-}(\vec r, \lvert \mk I \rvert)} \Eu^{\theta} (TM_{-}(\vec r, \lvert \mk I \rvert)) \\
\label{wc}
&=
 \sum_{\mk I^{\sharp} \in \mc D^{\ell}(\mk I)} \frac{ (n-d^{\sharp})! ( s( \mk I^{\sharp}, \mk I^{\flat}) - \bar{r}d^{\sharp} )  
 (\theta/\e+1)_{d^{\sharp}}}{ n! d^{\sharp}} \int_{\wt{M}^{\min(\mk I^{\sharp})-1} 
 (\vec r, \mk I^{\flat})} \tilde{\psi},
\end{align}
where $\bar{r} = r_{0} - r_{1}$ and we put $d^{\sharp}=\lvert \mk I^{\sharp} \rvert$ as before.

For $\mbi{\mk I}=( \mk I_{1}, \ldots, \mk I_{j}) \in \text{Dec}_{j}^{n}$, 
we use a symbol $\mbi{\mk I}_{> i} = \mk I_{\infty} \sqcup \bigsqcup_{k >i} \mk I_{k}$.
By substituting similar computations into \eqref{main>0} in Theorem \ref{thm:wcm}, we have
\begin{align}
&
\notag
\int_{M_{+}(\vec r, n)} \Eu^{\theta} (TM_{+}(\vec r, n))
-\int_{M_{-}(\vec r, n)} \Eu^{\theta} (TM_{-}(\vec r, n))  \\
\label{formula4}
&=
\sum_{j=1}^{n} \sum_{\mbi{\mk I} \in \Dec_{j}^{n}} \frac{ \lvert \mk I_{\infty} \rvert!  }{ n! }
\prod_{i=1}^{j} \frac{ ( s( \mk I_{i}, \mbi{\mk I}_{> i}) - \bar{r}\lvert \mk I_{i}\rvert )(\theta/\e +1)_{d_{i}}}{\lvert \mk I_{i}\rvert }
\int_{ M_{-}(\vec{r}, n - \lvert \mbi d_{\mbi{\mk I}} \rvert )} \psi.
\end{align}

\subsection{Vortex partition function}
\label{subsec:vortex}
In the following, we only consider the case where 
$J_{0} = \lbrace 1, \ldots, r_{0} \rbrace$ and 
$J_{1}=\lbrace r_{0}+1, \ldots, r \rbrace$.
Hence we use 
$Z^{\vec{r}}_{\pm \adj}$ instead of 
$Z^J_{\pm \adj}$.

We put $\alpha_{n} = \int_{M_{-}(\vec r, n)} \psi, \beta_{n} = \int_{M_{+}(\vec r, n)} \psi$, and rewrite 
\eqref{formula4} to analize vortex partition functions 
$Z^{\vec{r}}_{-\fund} (\e, \mbi a, \mbi m,  p) = \sum_{n=0}^{\infty} \alpha_{n} p^{n}$ and 
$Z^{\vec{r}}_{+\fund} (\e, \mbi a, \mbi m,  p) = \sum_{n=0}^{\infty} \beta_{n} p^{n}$.

We put $a_{\mbi{\mk I}}(\bar{r})=
\prod_{i=1}^{j} \frac{ (  s(\mk I, \mbi{\mk I}_{>i}) - \bar{r} \lvert \mk I_{i}\rvert) (\theta/\e+1)_{d_{i}}}{\lvert \mk I_{i}\rvert }$ 
for $\mbi{\mk I} =(\mk I_{1}, \ldots, \mk I_{j})$.
We also write by $\Dec (i, n)$ the set of elements 
$\mbi{\mk I}=(\mk I_{1}, \ldots, \mk I_{j})$ in $\Dec_{j}^{n}$ for $j = 1, \ldots, n$ such that $\lvert \mk I_{\infty} \rvert=i$.
For $i=0$, we put 
\[
A_{n} ( \bar{r})= \sum_{\mbi{\mk I} \in \Dec(0,n)} a_{\mbi{\mk I}}( \bar{r}).
\]
By computer experiments, we can check the following conjecture.
\begin{conj}
\label{conj}
For any $i=0, 1, \ldots, n$, we have 
$\sum_{\mbi{\mk I} \in \Dec(i, n)} a_{\mbi{\mk I}} (\bar{r})= \begin{pmatrix} n \\ i \end{pmatrix} A_{n-i}( \bar{r})$.
\end{conj}
In the following, we assume that this conjecture holds.
Substituting these into \eqref{formula4}, we get
\begin{align}
\beta_{n} - \alpha_{n}  
&=
\sum_{i=0}^{n} \frac{ i!  }{ n! } \sum_{\mbi{\mk I} \in \Dec(i,n)} 
\alpha_{\mbi{\mk I}} (\bar{r} )\alpha_{i}
\label{formula5}
=
\sum_{i=0}^{n} \frac{A_{n-i}( \bar{r})}{(n-i)!} \alpha_{i}.
\end{align}
This is equivalent to saying that
\begin{align}
\label{formula5'}
Z_{+}^{\vec{r}}(\e, \mbi a, \theta,  p ) = 
\left( \sum_{n=0}^{\infty} \frac{A_{n}(\bar{r})}{n!} p^{n} \right)
Z_{-}^{\vec{r}}(\e, \mbi a, \theta,  p). 
\end{align}

Since we have $Z_{+}^{(r, 0)}(\e, \mbi a, \theta,  p ) = Z_{-}^{(0, r)}(\e, \mbi a, \theta,  p)=1$,
from \eqref{formula5'} we get
\[
\sum_{n=0}^{\infty} \frac{A_{n}(\bar{r})}{n!} p^{n}
=
\begin{cases}
Z_{-}^{(\bar r, 0)}(\e, \mbi a, \theta,  p)^{-1} & \text{if } \bar r =r_{0}-r_{1}\ge 0, \\
Z_{+}^{(0, -\bar r)}(\e, \mbi a, \theta,  p) & \text{if } \bar r =r_{0} - r_{1} \le 0. 
\end{cases}
\]
In particular, we get $\sum_{n=0}^{\infty} \frac{A_{n}(\pm 1)}{n!} p^{n}
=(1-p)^{\pm (\theta/ \e +1)}$, and
\begin{align}
\label{formula6}
 Z_{-}^{(r_{0}, r_{0} \pm 1)}(\e, \mbi a, \theta,  p)
=
(1-p)^{\pm (\theta/\e+1)}
Z_{+}^{(r_{0}, r_{0} \pm 1)}(\e, \mbi a, \theta,  p ) 
\end{align}

We take limit $a_{r}/\e \to \infty$ in \eqref{comb3} and \eqref{comb4}:
\begin{align*}
\lim_{a_{r}/\e \to \infty} 
Z^{(r_{0}, r_{1})}_{-}  (\e, \mbi a, \theta, p )
&=
Z^{(r_{0}, r_{1}-1)}_{-}  (\e, a_{1}, \ldots, a_{r-1}, \theta, p )\\
\lim_{a_{r} / \e \to \infty} 
Z^{(r_{0}, r_{1})}_{+} (\e, \mbi a, \theta, p ) 
&=
(1- p)^{-\theta/\e - 1}
Z^{(r_{0}, r_{1}-1)}_{+} (\e, a_{1}, \ldots, a_{r-1}, \theta, p ).
\end{align*}
Applying this to \eqref{formula6}, we have 
$Z^{(r_{0}, r_{0})}_{-}  (\e, \mbi a, \theta, p )=
Z^{(r_{0}, r_{0})}_{+} (\e, \mbi a, \theta, p )$, and
\[
Z^{(r_{0}, r_{0}-2)}_{-}  (\e, a_{1}, \ldots, a_{r-1}, \theta, p )=
(1- p)^{-2(\theta/\e+1)}
Z^{(r_{0}, r_{0}-2)}_{+} (\e, a_{1}, \ldots, a_{r-1}, \theta, p ) .
\]
Repeating this procedure, we have
\[
Z^{(r_{0}, r_{0}-i)}_{-} (\e, \mbi a, \theta, p ) 
=
(1- p)^{-i(\theta/\e+1)}
Z^{(r_{0}, r_{0}-i)}_{+} (\e, \mbi a, \theta, p )
\] 
for any positive integer $i \le r_{0}$. 
In particular, we have $Z^{(r_{0}, 0)}_{-} (\e, \mbi a, \theta, p )
=(1- p)^{-r_{0} (\theta / \e +1) }$.
Similarly we have
\[
Z^{(0, r_{1})}_{+} (\e, \mbi a, \theta, p )= 
\sum_{n=0}^{\infty} \frac{A_{n}(-r_{1})}{n!} p^{n} = 
(1- p)^{-r_{1} (\theta/\e+1)}.
\]
This give the following proposition.
\begin{rem}
For any integer $\bar{r}$, we have 
\[
\sum_{n=0}^{\infty} \frac{A_{n}(\bar{r})}{n!} p^{n}
=
(1-p)^{\bar{r} (\theta/ \e +1)}.
\]
\end{rem}
Thus we get Theorem \ref{main} (b).

\begin{NB}
\section{Maple check}
\subsection{Check of Conjecture \ref{conj}}
We want to check that for any $0 \le k \le n$, we have 
\begin{align*}
&
\sum_{j=1}^{k} \sum_{\substack{\mbi{\mk I} \in \Dec_{j}^{n} \\ |\mbi d_{\mbi{\mk I}}|=k }} 
\prod_{l=1}^{j} \frac{ ( s( \mk I_{l}, \mbi{\mk I}_{> l}) - r \lvert \mk I_{l}\rvert )(u)_{d_{l}}}{\lvert \mk I_{l}\rvert }
\\
=&
\binom{n}{k}
\sum_{j=1}^{k} \sum_{\substack{\mbi{\mk I} \in \Dec_{j}^{k} \\ |\mbi d_{\mbi{\mk I}}|=k}} 
\prod_{l=1}^{j} \frac{ ( s( \mk I_{l}, \mbi{\mk I}_{> l}) - r \lvert \mk I_{l}\rvert )(u)_{d_{l}}}{\lvert \mk I_{l}\rvert }
\\
=&
(-ru)_{k}.
\end{align*}

\begin{verbatim}



poch:=proc(x, k)
local f, i;    
f:=1; 
for i from 0 to k-1do 
f:=f*(x +i) od;
RETURN(f);
end:



EWop:=(y,m)->[op(y),m]:



ListSetPtns:=proc(n,k)
local east,west,i,out ; options remember;
#Lists all partitions of the set 1,..,n into k classes 
#Output array[i] is the class to which i belongs.
if n=1 then
if k<>1 then RETURN([]) else RETURN([[1]]) fi: else
east:=ListSetPtns(n-1,k-1): west:=ListSetPtns(n-1,k): out:=map(EWop,east,k);
for i from 1 to k do out:=[op(out),op(map(EWop,west,i))] od; 
RETURN(out);
fi: end:


dec:=proc(n,j)
local i, k, L, M, N;
N:=[];
for L in ListSetPtns(n, j+1) do 
for k from 1 to j+1 do
M := [seq([], i = 1 .. j+1)];
for i from 1 to n do
if L[i] < k then
M[L[i]]:=[op(M[L[i]]), i];
fi;
if L[i] > k then
M[L[i]-1]:=[op(M[L[i]-1]), i];
fi;
if L[i] = k then
M[j+1]:=[op(M[j+1]), i];
fi;
od;
M:=[seq(M[j-i+1], i=1..j), M[j+1]];
N:=[op(N), M];
od; od;
for L in ListSetPtns(n, j) do 
M := [seq([], i = 1 .. j+1)];
for i from 1 to n do
M[L[i]]:=[op(M[L[i]]), i];
od;
M:=[seq(M[j-i+1], i=1..j), M[j+1]];
N:=[op(N), M];
od;
return N;
end;

dimv:=proc(L)
return map(l->nops(l), L);
end;


dim := proc(L) 
local d, i, f; 
d:=dimv(L);
f := 0; 
for i to nops(d) - 1 do f := f + d[i]; end do; 
return f; 
end;


s := proc(L, l) 
local i, j, k, x, y; 
x := 0; y := 0;
for i in L[l] do  
for j from 1 to nops(L) - l do 
for k in L[j+l] do 
if i < k  then x := x + 1; end if; 
if i > k then y := y + 1; end if;  
end do; end do; end do; 
RETURN(x-y); 
end;


check1 := proc(d, r) 
local i, j, f, g, l, n, L; 
f := 0; 
n := sum(d[i], i = 1 .. nops(d)); 
j := nops(d) - 1; 
for L in dec(n, j) do if dimv(L) = d then 
g := 1; 
for l to j do 
g := g*(s(L, l) - r*nops(L[l])); 
end do; 
f := f + g; 
end if; end do; 
return factor(f); 
end;


check2:=proc(k, n)
local  i,j, f, g, l, L;
f := 0;
for j from 1 to k do
for L in dec(n, j) do if dim(L) = k then 
g:=1;
for l from 1 to j do
g := g* (s(L, l) - nops(L[l]))*poch(u, nops(L[l]))/nops(L[l]); 
end do;
f:=f+g;
end if; end do;
end do;
return(factor(f- binomial(n,k)*poch(-u, k)));
end;



\end{verbatim}

\subsection{Check of \eqref{check}}
To prove \eqref{adj} for general $\bar{r}$, we want to check
\begin{align}
\label{check}
\sum_{\mbi{\mk I} \in \Dec(0,n)} \prod_{k=1}^{j} \frac{ ( s(k, \mbi{\mk I}) - t |\mk I_{k}| ) (u)_{d_{k}}}{|\mk I_{k}| }
=(- t u)_{n},
\end{align}
where $s(k, \mbi{\mk I}) =s( \mk I_{k}, \mk I_{k+1} \sqcup \cdots \sqcup \mk I_{j} )$, and 
\[
s(\mk I_{1}, \mk I_{2}) = \left| \lbrace (i, j) \in \mk I_{1} \times \mk I_{2} \mid i < j \rbrace \right| - \left| \lbrace (i, j) \in \mk I_{1} \times \mk I_{2} \mid i > j \rbrace \right|.
\]

This implies $A_{n} (\bar{r} )= (- \bar{r}u)_{n}$, and  we get the desired equality.
\[
Z_{+}^{J}(\e, \mbi a, \theta,  p)  = (1-p)^{\bar{r} u} Z_{-}^{J}(\e, \mbi a, \theta,  p).
\]

We want to check
\[
\sum_{n=0}^{\infty} \frac{A_{n}(\bar{r})}{n!} p^{n}
=
(1-p)^{\bar{r} u}
\]
without assumption $|\bar{r}|=|r_{0}-r_{1}|\le 1$.

We recall $u=(\theta+\e)/\e$ and 
we also write by $\Dec (0, n)$ the set of tuples $\mk I^{\sharp}=(\mk I_{1}, \ldots, \mk I_{j})$ of non-empty subsets of $[n]$ such that $\mk I_{1} \sqcup \cdots \sqcup \mk I_{j}=[n]$ and $\min(\mk I_{1}) > \cdots > \min(\mk I_{j})$.
We put
\[
A_{n} ( \bar{r})= \sum_{\mk I^{\sharp} \in \Dec(0,n)} \prod_{k=1}^{j} \frac{ ( s( \mk I_{k}, \mk I_{k+1} \sqcup \cdots \sqcup \mk I_{j} ) - \bar{r} |\mk I_{k}| ) (u)_{d_{k}}}{|\mk I_{k}| },
\]
where $s(\mk I_{1}, \mk I_{2}) = \left| \lbrace (i, j) \in \mk I_{1} \times \mk I_{2} \mid i < j \rbrace \right| - \left| \lbrace (i, j) \in \mk I_{1} \times \mk I_{2} \mid i > j \rbrace \right|$.

\begin{verbatim}



EWop:=(y,m)->[op(y),m]:



ListSetPtns:=proc(n,k)
local east,west,i,out ; options remember;
#Lists all partitions of the set 1,..,n into k classes #Output array[i] is the class to which i belongs.
if n=1 then
if k<>1 then RETURN([]) else RETURN([[1]]) fi: else
east:=ListSetPtns(n-1,k-1): west:=ListSetPtns(n-1,k): out:=map(EWop,east,k);
for i from 1 to k do out:=[op(out),op(map(EWop,west,i))] od; 
RETURN(out);
fi: end:



dimv := proc(L, j) 
local i, M; 
M := [seq(0, i = 1 .. j )]; 
for i from 1 to nops(L) do 
if L[i] <= j then
M := applyop(x -> x + 1, L[i], M); 
fi;
end do; 
return M; 
end;



s := proc(L, l) 
local i, j, x, y; 
x := 0; y := 0; 
for i from 1 to nops(L) - 1 do for j from i + 1 to nops(L) do 
if L[i] = l then if l < L[j] or L[j]=0 then x := x + 1; end if; end if; 
if L[j] = l then if l < L[i] or L[j]=0 then y := y + 1; end if; end if; 
end do; end do; 
RETURN(x-y); 
end;



check:=proc(d, n, r)
local  i,j, f, k, l, m, L, M;
# d=(d_{1}, \ldots, d_{j})
j:=nops(d);
f:=0;
for L in ListSetPtns(n, j+1) do 
for k from 1 to j+1 do
M:=[];
for i from 1 to n do
if L[i] < k then
M:=[op(M), L[i]];
fi;
if L[i] > k then
M:=[op(M), L[i]-1];
fi;
if L[i] = k then
M:=[op(M), j+1];
fi;
od;
if dimv(M,j)= d then
M:=[seq(j+1-M[i], i=1..n)];
f:=f+product(s(M, l), l=1..j);
fi;
od;
od;
for L in ListSetPtns(n, j) do 
if dimv(L,j)= d then
L:=[seq(j+1-L[i], i=1..n)];
f:=f+product(s(L, l) - r * d[l], l=1..j);
fi;
od;
RETURN(factor(f));
end;




\end{verbatim}

\end{NB}


\noindent {\bf Acknowledgements}\\
The authors thanks Hidetoshi Awata, Ayumu Hoshino, Hiroaki Kanno, Hitoshi Konno, Takuro Mochizuki, Hiraku Nakajima, Masatoshi Noumi, Takuya Okuda,  
Yusuke Ohkubo, Yoshihisa Saito and Junichi Shiraishi for discussion.
He is grateful for Masatoshi Noumi for directing his attention to relationships between our results and transformation formulas for
multiple hypergeometric functions including the Kajihara transformations.
RO is partially supported by Grant-in-Aid for Scientific Research 21K03180 and 17H06127, JSPS.
He has had the generous support and encouragement of Masa-Hiko Saito.
YY is supprted by Grant-in-Aid for Scientific Research 21K03382 and 21H05190, JSPS.
This work was partly supported by Osaka Central Advanced Mathematical
Institute: MEXT Joint Usage/Research Center on Mathematics and
Theoretical Physics JPMXP0619217849, and by the Research Institute for Mathematical Sciences,
an International Joint Usage/Research Center located in Kyoto University.
The authors have no conflicts of interest directly relevant to the content of this article.
Data sharing not applicable to this article as no datasets were generated or analysed during the current study.

\end{document}